\documentclass[11pt]{article}
\usepackage{amsmath}
\usepackage{amsthm}
\usepackage{fullpage}
\usepackage{graphicx}
\usepackage{pstricks}
\newpsobject{showgrid}{psgrid}{subgriddiv=1,griddots=10,gridlabels=6pt}

\newcommand {\abs} [1] {\left| #1 \right|}
\newcommand {\Abs} [1] {\bigl\lvert #1 \bigr\rvert}
\newcommand {\ABS} [1] {\Biggl\lvert #1 \Biggr\rvert}
\newcommand {\ie} {i.e.}
\newcommand {\E}{{\mathrm{E}}}
\renewcommand{\P}{{\mathrm{P}}}
\newcommand{\var}{\mathop{\mathrm{Var}}}
\newcommand{\covar}{\mathop{\mathrm{Covar}}}
\newcommand{\des}{\mathop{\mathrm{des}}}

\newcommand{\nghd}[1]{\mathop{\mathrm{nghd}}\{#1\}}

\newcommand {\euler}[2]{\left\langle\begin{smallmatrix} #1\\#2
                        \end{smallmatrix}\right\rangle}
\newcommand {\Euler}[2]{\left<\begin{array}{c} #1\\#2
                        \end{array}\right>}
\newtheorem {thm} {Theorem}[section]

\newtheorem {lem} [thm] {Lemma}

\theoremstyle {definition}

\theoremstyle {remark}

\numberwithin{equation}{section}
\title{Shuffling Cards for Blackjack, Bridge, and Other Card Games}
\author{Mark Conger and D. Viswanath \thanks{ Department of Mathematics, University of Michigan, 
530 Church Street, Ann Arbor, MI 48109, U.S.A.,
{\tt mconger@umich.edu and divakar@umich.edu}.
This work was supported by a research fellowship from the Sloan Foundation.
Computations on the Teragrid networks at NCSA and SDSC were supported by
a DAC award.}}
\date{March 23, 2006}
\begin{document}
\maketitle

\begin{abstract}
This paper is about the following question: How many riffle shuffles
mix a deck of card for games such as blackjack and bridge?  An object
that comes up in answering this question is the descent polynomial
associated with pairs of decks, where the decks are allowed to have
repeated cards. We prove that the problem of computing the descent
polynomial given a pair of decks is $\#P$-complete. We also prove that
the coefficients of these polynomials can be approximated using the
bell curve. However, as must be expected in view of the
$\#P$-completeness result, approximations using the bell curve are not
good enough to answer our question. Some of our answers to
the main question are supported
by theorems, and others are based on experiments supported by
heuristic arguments.  In the introduction, we carefully discuss the
validity of our answers.
\end{abstract}

\section{Introduction}
\subsection{Probability theory, card shuffling, and computational complexity theory}
The representation of a sequence of coin tosses as the binary digits
of a real number appears in a 1909 paper by Borel \cite{Borel1}. This
representation was a step towards the formalization of the notion of
probability using measure theory. Card shuffling was another example
that was discussed around that time by Borel, Poincar\'{e}, and
others. However, the analysis of card shuffling, unlike that of coin
tosses, was not advanced very far at that time.  Even as he
noted that a large number of shuffles brings the distribution of a
deck close to the uniform distribution, R\'{e}nyi wrote in the 1960s
that he would not deal with ``what is meant by a large enough number
of movements'' \cite{RenyiBook}. Much is now known about the number of
shuffles for mixing a deck of cards thanks to the development of the
convergence theory of Markov chains by Aldous, Diaconis, and others
\cite{Aldous1} \cite{DiaconisBook}.

A notable result about riffle shuffles is due to Bayer and Diaconis
\cite{BD1}. The riffle shuffle, where a deck of cards is cut into two
packets and cards are dropped from the two packets in some order, is
the most common method of shuffling cards. For a certain model of the riffle
shuffle, Bayer and Diaconis gave a complete analysis of 
riffle shuffles. This result assumes that all $52$ cards are
distinct and that each of the $52!$ possible permutations must be
nearly equally likely for the deck to considered well mixed.

However, in card games such as blackjack, the distinction between
suits (a standard deck of $52$ cards has $4$ suits with cards labeled
$\mathtt{A, 2,\ldots, 10, J, Q, K}$ in each suit) is ignored. In other card
games such as bridge, all $52$ cards are distinct, but there are only
$4$ players. In bridge, it is enough if each player receives a random
{\it set} of $13$ cards --- the order in which those $13$ cards are received
by a player is inconsequential. Therefore, in some card games, not all
cards in the standard deck are distinct, and in others, not all $52!$
permutations need to be nearly equally likely for the deck to be
considered well mixed. The mixing times for such card games is the
topic of this paper.

The popularity of card games is a reason to study card shuffling, but
it is not the only reason. Examples based on card shuffling are an
inextricable part of the convergence theory of finite state Markov
chains. One application of this theory is to theoretical computer
science, and in particular, to the problem of computing the permanent
of a matrix all of whose entries are $0$ or $1$.

 Let $A$ be an $n\times n$ matrix whose $ij$th entry is
$a_{ij}$. Then its permanent is defined as
\begin{equation}
\mathop{\mathrm{perm}}(A) = \sum_\pi a_{1\pi(1)} a_{2\pi(2)}\ldots a_{n\pi(n)},
\label{eqn-1-1}
\end{equation}
where $\pi$ ranges over the permutations of $\{1,2,\ldots,n\}$. If the
term under the summation in \eqref{eqn-1-1} were multiplied by $\pm
1$, with the sign being $+$ for even $\pi$ and $-$ for odd $\pi$, we
would have a definition of the determinant of $A$.  In spite of this
resemblance to the determinant, computing the permanent is much harder
than computing the determinant. As every student of linear algebra
knows, the determinant can be computed using $O(n^3)$ arithmetic
operations by carrying out Gaussian elimination or row reduction. But
there is no known way of computing the permanent that is much better
than using its definition \eqref{eqn-1-1} directly. The operation
count for the direct method is more than $n!$. Even for $n=50$, such
an operation count is far out of the reach of today's computers.

The permanent belongs to a class of counting problems known as
$\#P$. It is also $\#P$-complete which means that every problem in
$\#P$ can be reduced to it in time that is polynomial in the size of
the problem. There are many other $\#P$-complete problems, some of
them of considerable importance to trade and industry, but it is
believed that there is no polynomial time algorithm for any of these
problems.  The $P \neq NP$ conjecture is closely related to this
belief.  Proof or disproof of either conjecture would be a major advance in
computational complexity theory.

If it is only desired to approximate the permanent, rather than
compute it exactly, Markov chains are of help. Jerrum, Sinclair, and
Vigoda \cite{JSV1} have devised a Markov chain, analyzed its
convergence, and proved that the permanent can be approximated
accurately with high probability in polynomial time.

This relationship between Markov chains and computational complexity
theory is inverted in our work. To find out the number of riffle
shuffles that mix a deck of cards for blackjack, bridge and other card
games, we need to be able to find the transition probability between
two given decks, with some cards repeated, under a given number of
riffle shuffles.  There are simple formulas to pass back and forth
between transition probabilities between two decks and a polynomial
that will be called the descent polynomial for those two decks.  We
prove that given two decks, the problem of computing their descent
polynomial is $\#P$-complete.

A graph of the coefficients of the descent polynomial looks strikingly
like a bell curve for most pairs of decks. We prove a theorem
showing that coefficients of the descent polynomial are approximated
by the normal law in very general circumstances. The coefficients of
the descent polynomial that we need to figure out the mixing times for
games such as blackjack and bridge lie in the tail of the normal fit.
Unfortunately, but not unexpectedly, the bounds on the normal
approximation are not sharp enough in this region. The proof of 
$\#P$-completeness in Section 4 suggests that these must be the
coefficients that are hard to compute accurately, and indeed they are.
It would be too simplistic to expect a probablistic technique useful
for proving normal approximation to provide a way around a 
$\#P$-complete problem.

We find a way around using computations justified by heuristic
arguments.  The normal approximation result helps in two ways --- it
suggests a probablistic approach to finding the descent polynomial
given a pair of decks and our computations of the descent polynomial
are guided by the form of the bell curve. We state some of our results and
discuss the validity of our computations in the second part of this
introduction.

Apart from computational complexity theory and probability theory, the
analysis of the mixing time for games such as blackjack and bridge is
also connected to the theory of descents. The definition of descents
of permutations is given in the next section and the central role of
descents in the analysis of riffle shuffles will become clear. The
systematic study of descents was begun by MacMahon
\cite{MacMahonBook}.  More recent references on this topic are Foata
and Sch\"{u}tzenberger \cite{FSBook}, Gessel and Reutenauer
\cite{GeRe1}, and Knuth \cite{KnuthVol3}.

\subsection{Theorems, experiments, and card games}

\begin{table}
\begin{center}
\scriptsize
\begin{tabular}{c|c|c|c|c|c|c|c|c|c|c}
 & 1 & 2 & 3 & 4 & 5 &  6 & 7 & 8 & 9 & 10\\ \hline
&&&&&&&&&&\\
BayerDiaconis &  {\bf 1} & {\bf 1} & {\bf 1} & {\bf 1} & 
{\bf .924} & {\bf .614} & {\bf .334} & {\bf .167} & {\bf .085} & {\bf .043}\\ 
&&&&&&&&&&\\
\hline
&&&&&&&&&&\\
Blackjack1 & {\bf 1} & {\bf 1} & {\bf 1} & {\bf .481} & 
{\bf .215}  & {\bf .105}  & {\bf .052}  & {\bf .026}  & {\bf .013}  & {\bf .007}\\ 
&&&&&&&&&&\\
\hline
&&&&&&&&&&\\
Bridge2 & & & .45? & .16 & 
.08 & .04 & .02 & .01 & .00 & .00\\ 
&&&&&&&&&&\\
\end{tabular}
\normalsize
\end{center}
\caption[xyz]{Table of total variation distances after $1$ to $10$ riffle
shuffles for three scenarios. The meaning and validity of these numbers
is discussed in the text.}
\label{table-0}
\end{table}

Each of the three lines of Table \ref{table-0} corresponds to a
different scenario. In the first scenario, which was considered by
Bayer and Diaconis \cite{BD1}, there are $52$ distinct cards and we
want all $52!$ permutations to be nearly equally likely. In the second
scenario (blackjack), the distinction between the suits is ignored and
the source deck is assumed to be four aces on top of four $2$s, on
top of four $3$s, and so on. We want all $52!/4!^{13}$ possible
permutations to be nearly equally likely. In the third scenario
(bridge), all $52$ cards are distinct but they are dealt to four
players in cyclic order, and we only want the partition of the
$52$ cards to the four players to be random.

The total variation distances in Table \ref{table-0} are a measure
of how well mixed the deck is after a certain number of riffle shuffles.
After $7$ riffle shuffles the total variation distance is $0.334$
for the Bayer-Diaconis scenario. The total variation distance is
lesser after only $5$ riffle shuffles for blackjack and after
only $4$ riffle shuffles for bridge.

As already mentioned, we had to resort to careful computations and
heuristic arguments to determine some of the numbers in Table
\ref{table-0}.  How reliable are those numbers? How rigorous are the
methods used to find them?   We now answer these questions.

Each number in the first line of Table \ref{table-0} is given by
a formula derived by Bayer and Diaconis \cite{BD1}. The formula
is quite simple to implement and the rounding errors in finding
the numbers in the first line can be bounded easily. We are happy
to accept the numbers in that line as theorems.

The numbers in the second line of Table \ref{table-0} are about
blackjack. We prove that each blackjack number has an error
less than $.001$ with a probability greater than $96\%$, or an
error less than $.01$ with a probability greater than $99.9996\%$.

The proofs of those error estimates are valid only with an assumption,
however. The numbers are generated after permuting the blackjack deck
randomly $10$ million times, and the proofs of the error estimates
assume those permutations to be independent of each other. The issue
is whether the computer generated pseudorandom numbers match the
idealizations found in probability theory.

Pseudorandom numbers have been studied and used extensively
\cite{KnuthVol2}. A specific example is given by the coupling from the
past construction by Propp and Wilson
\cite{PW}. They prove that a sequence of random steps taken in a
particular way is guaranteed to terminate at a state of the Ising
model with a particular probability distribution. Implementations of
that method, however, make do with pseudorandom numbers. The situation
with our error estimates is similar.

Finally, what of the bridge numbers found in the third line of
Table \ref{table-0}? We prove no error estimates for these numbers.
While the blackjack problem fits into a subclass of pairs of
decks for which we can find the descent polynomial quickly
\cite{CV1}, the bridge problem does not. Considering that
the problem of finding the descent polynomial given any pair of decks
is $\#P$-complete and that ``counting problems that can be solved
exactly in polynomial time are few and far between'' (see
\cite{Jerrum}), one begins to suspect that theorems that assert error
bounds for the bridge numbers might be out of reach for a long time.

The bridge numbers in Table \ref{table-0} were obtained using careful
computations guided by heuristic arguments. These numbers
have the same type of validity as the numbers produced in experimental
physics and chemistry. The method for producing the bridge numbers
has been checked for internal consistency in a variety of ways in
Section 7, and the method is reported in enough detail to permit
others to reproduce our results. The numbers are open to refutation,
a quality of experimental results that has sometimes been emphasized
\cite{Popper}.

\section{A model of riffle shuffles}
The figure below shows a riffle shuffle of a deck of $5$ cards
numbered $1$ through $5$.
\begin{center}
\psset{xunit=1.200cm,yunit=0.400cm}
\begin{pspicture}(10.000,6.000)
\rput*(0.000,5.400){1}
\psframe[fillstyle=solid,fillcolor=white](0.200,5.320)(1.600,5.480)
\rput*(0.000,4.200){2}
\psframe[fillstyle=solid,fillcolor=white](0.200,4.120)(1.600,4.280)
\rput*(0.000,3.000){3}
\psframe[fillstyle=solid,fillcolor=white](0.200,2.920)(1.600,3.080)
\rput*(0.000,1.800){4}
\psframe[fillstyle=solid,fillcolor=black](0.200,1.720)(1.600,1.880)
\rput*(0.000,0.600){5}
\psframe[fillstyle=solid,fillcolor=black](0.200,0.520)(1.600,0.680)
\psline[linestyle=dashed](-0.200,2.400)(1.900,2.400)
\rput(2.45,4.5){Cut}
\psline[linewidth=.1cm]{->}(2.2,3)(2.7,3)
\rput*(2.950,4.200){1}
\psframe[fillstyle=solid,fillcolor=white](3.150,4.120)(4.550,4.280)
\rput*(2.950,3.000){2}
\psframe[fillstyle=solid,fillcolor=white](3.150,2.920)(4.550,3.080)
\rput*(2.950,0.600){3}
\psframe[fillstyle=solid,fillcolor=white](3.150,0.520)(4.550,0.680)
\rput*(4.800,5.400){4}
\psframe[fillstyle=solid,fillcolor=black](5.000,5.320)(6.400,5.480)
\rput*(4.800,1.800){5}
\psframe[fillstyle=solid,fillcolor=black](5.000,1.720)(6.400,1.880)
\rput(6.95,4.5){Riffle}
\psline[linewidth=0.1cm]{->}(6.7,3)(7.2,3)
\rput*(7.600,5.400){4}
\psframe[fillstyle=solid,fillcolor=black](7.800,5.320)(9.200,5.480)
\rput(9.500,5.400){(2)}
\rput*(7.600,4.200){1}
\psframe[fillstyle=solid,fillcolor=white](7.800,4.120)(9.200,4.280)
\rput(9.500,4.200){(1)}
\rput*(7.600,3.000){2}
\psframe[fillstyle=solid,fillcolor=white](7.800,2.920)(9.200,3.080)
\rput(9.500,3.000){(1)}
\rput*(7.600,1.800){5}
\psframe[fillstyle=solid,fillcolor=black](7.800,1.720)(9.200,1.880)
\rput(9.500,1.800){(2)}
\rput*(7.600,0.600){3}
\psframe[fillstyle=solid,fillcolor=white](7.800,0.520)(9.200,0.680)
\rput(9.500,0.600){(1)}
\end{pspicture}

\end{center}
The first step in a riffle shuffle is the {\it cut}. In the picture
above, the deck is cut below the third card to get two packets with
$3$ and $2$ cards, respectively. The second step is the {\it
riffle}. The packets are riffled by repeatedly 
dropping cards from the bottom of
one of the two packets until a new shuffled deck is obtained.  In the
picture above, the second card to be dropped and the
last card to be dropped are from the second packet and all
others are from the first, as indicated by the
parenthesized numbers on the right.

In a riffle shuffle, the cut can be placed anywhere in the deck and the
two packets can be riffled together in different ways. To model riffle
shuffles, it is necessary to assign probabilities to the various ways
of riffle shuffling a deck of $n$ cards. 

In the model we use, a random riffle shuffle of a deck of $n$ cards is
obtained by first generating a sequence of $n$ independent random
numbers each of which is either $1$ or $2$ with probability $1/2$.
For the riffle shuffle depicted above, the corresponding sequence is
given in parenthesis on the right.
The
cut must be placed such that the number of cards in the first packet
equals the number of $1$s in the sequence, and the number in the
second packet equals the number of $2$s.  During the riffle, a card
must be dropped from the first or second packet to get the $i$th card
in the shuffled deck according as the $i$th number in the sequence is $1$
or $2$. The riffle shuffle depicted above results from the sequence
$2,1,1,2,1$.

\begin{figure}
\begin{center}
\psset{xunit=1.000cm,yunit=0.600cm}
\begin{pspicture}(12.000,6.000)
\rput*(0.000,5.625){1}
\psframe[fillstyle=solid,fillcolor=white](0.200,5.545)(1.600,5.705)
\rput*(0.000,4.875){2}
\psframe[fillstyle=solid,fillcolor=white](0.200,4.795)(1.600,4.955)
\rput*(0.000,4.125){3}
\psframe[fillstyle=solid,fillcolor=gray](0.200,4.045)(1.600,4.205)
\rput*(0.000,3.375){4}
\psframe[fillstyle=solid,fillcolor=gray](0.200,3.295)(1.600,3.455)
\rput*(0.000,2.625){5}
\psframe[fillstyle=solid,fillcolor=gray](0.200,2.545)(1.600,2.705)
\rput*(0.000,1.875){6}
\psframe[fillstyle=solid,fillcolor=black](0.200,1.795)(1.600,1.955)
\rput*(0.000,1.125){7}
\psframe[fillstyle=solid,fillcolor=black](0.200,1.045)(1.600,1.205)
\rput*(0.000,0.375){8}
\psframe[fillstyle=solid,fillcolor=black](0.200,0.295)(1.600,0.455)
\psline[linestyle=dashed](-0.200,4.500)(1.900,4.500)
\psline[linestyle=dashed](-0.200,2.250)(1.900,2.250)
\rput(2.45,3.6){Cut}
\psline[linewidth=.1cm]{->}(2.2,3)(2.7,3)
\rput*(2.950,5.625){1}
\psframe[fillstyle=solid,fillcolor=white](3.150,5.545)(4.550,5.705)
\rput*(2.950,0.375){2}
\psframe[fillstyle=solid,fillcolor=white](3.150,0.295)(4.550,0.455)
\rput*(4.800,4.125){3}
\psframe[fillstyle=solid,fillcolor=gray](5.000,4.045)(6.400,4.205)
\rput*(4.800,2.625){4}
\psframe[fillstyle=solid,fillcolor=gray](5.000,2.545)(6.400,2.705)
\rput*(4.800,1.125){5}
\psframe[fillstyle=solid,fillcolor=gray](5.000,1.045)(6.400,1.205)
\rput*(6.650,4.875){6}
\psframe[fillstyle=solid,fillcolor=black](6.850,4.795)(8.250,4.955)
\rput*(6.650,3.375){7}
\psframe[fillstyle=solid,fillcolor=black](6.850,3.295)(8.250,3.455)
\rput*(6.650,1.875){8}
\psframe[fillstyle=solid,fillcolor=black](6.850,1.795)(8.250,1.955)
\rput(8.95,3.6){Riffle}
\psline[linewidth=0.1cm]{->}(8.7,3)(9.2,3)
\rput*(9.600,5.625){1}
\psframe[fillstyle=solid,fillcolor=white](9.800,5.545)(11.200,5.705)
\rput(11.500,5.625){(1)}
\rput*(9.600,4.875){6}
\psframe[fillstyle=solid,fillcolor=black](9.800,4.795)(11.200,4.955)
\rput(11.500,4.875){(3)}
\rput*(9.600,4.125){3}
\psframe[fillstyle=solid,fillcolor=gray](9.800,4.045)(11.200,4.205)
\rput(11.500,4.125){(2)}
\rput*(9.600,3.375){7}
\psframe[fillstyle=solid,fillcolor=black](9.800,3.295)(11.200,3.455)
\rput(11.500,3.375){(3)}
\rput*(9.600,2.625){4}
\psframe[fillstyle=solid,fillcolor=gray](9.800,2.545)(11.200,2.705)
\rput(11.500,2.625){(2)}
\rput*(9.600,1.875){8}
\psframe[fillstyle=solid,fillcolor=black](9.800,1.795)(11.200,1.955)
\rput(11.500,1.875){(3)}
\rput*(9.600,1.125){5}
\psframe[fillstyle=solid,fillcolor=gray](9.800,1.045)(11.200,1.205)
\rput(11.500,1.125){(2)}
\rput*(9.600,0.375){2}
\psframe[fillstyle=solid,fillcolor=white](9.800,0.295)(11.200,0.455)
\rput(11.500,0.375){(1)}
\end{pspicture}
\end{center}

\caption[xyz]{A $3$-shuffle of $8$ cards. The numbers in parentheses near the
right edge of the figure give the packet the card comes from.}
\label{fig-1}
\end{figure}
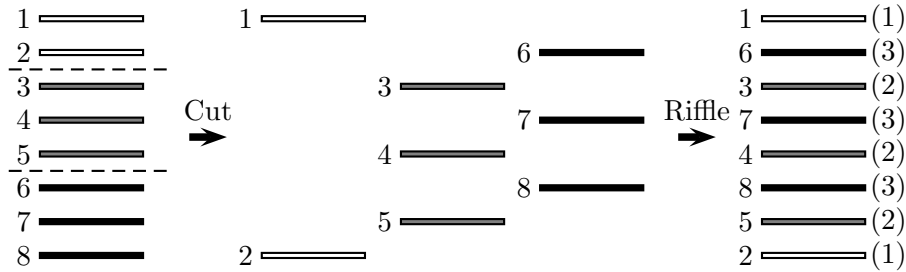

This model can be extended to assign probabilities to {\it
$a$-shuffles}, in which the deck is cut into $a$ packets, for any positive
integer $a$.  Figure \ref{fig-1} depicts a $3$-shuffle of $8$ cards. 
The extension is as follows. First, consider a sequence of $n$ random
numbers each of which is uniformly distributed over the set $\{1,
2,\ldots, a\}$ and independent of the $n-1$ other random numbers. When
the deck of $n$ cards is cut into $a$ packets, the number of cards in
the $p$th packet must be equal to the number of $p$s in the random
sequence.  If the number in the $i$th
position of the random sequence is $p$, then the card that ends up in
that position in the shuffled deck must be dropped from the $p$th
packet.  The random sequence that corresponds to the $3$-shuffle
depicted in Figure \ref{fig-1} is $1,3,2,3,2,3,2,1$. From here onwards,
the phrases {\it riffle shuffle} and {\it 2-shuffle} are used 
interchangeably.

In this model of the $a$-shuffle, the $p$th packet will be empty if
$p$ does not appear at all in the random sequence. If every instance
of $p$ precedes every instance of $p+1$ in the random sequence, then
all cards from the $p+1$st packet must be dropped before any card is
dropped from the $p$th packet. An $a$-shuffle has $a-1$ cuts, and
when some packets are empty, some of these cuts fall between the same
positions.

According to Bayer and Diaconis \cite{BD1}, this model was described
by Gilbert and Shannon in 1955 and independently by Reeds in 1971. It
is also described by Epstein \cite{EpsteinBook}, who assumes 
$a=2$ and calls it the
amateur shuffle.  It is sometimes called the GSR-model, after
three of the people who proposed it.
Many aspects of this model became clear only with the
work of Bayer and Diaconis.

One of the descriptions of this model of the $2$-shuffle given by
Bayer and Diaconis \cite{BD1} is as follows. A deck of $n$ cards
is cut after the $k$th card with probability $\frac{1}{2^n} \binom{n}{k}$;
in other words, the cut is binomial. During the riffle, if the first
packet has $a$ cards and the second packet has $b$ cards, the next
card to be dropped is from the bottom of the first packet with probability
$a/(a+b)$, and from the bottom of the second packet with probability
$b/(a+b)$. 

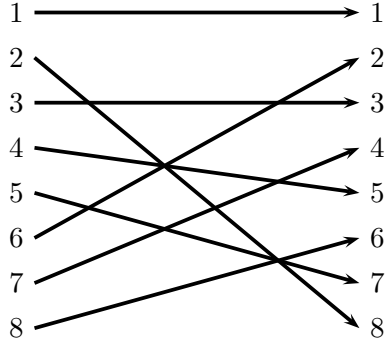
\begin{figure}
\begin{center}
\psset{xunit=1.200cm,yunit=0.400cm}
\begin{pspicture}(4.000,12.000)
\rput*(0.000,11.250){1}
\rput*(4.000,11.250){1}
\rput*(0.000,9.750){2}
\rput*(4.000,9.750){2}
\rput*(0.000,8.250){3}
\rput*(4.000,8.250){3}
\rput*(0.000,6.750){4}
\rput*(4.000,6.750){4}
\rput*(0.000,5.250){5}
\rput*(4.000,5.250){5}
\rput*(0.000,3.750){6}
\rput*(4.000,3.750){6}
\rput*(0.000,2.250){7}
\rput*(4.000,2.250){7}
\rput*(0.000,0.750){8}
\rput*(4.000,0.750){8}
\psline[linewidth=0.5mm]{->}(0.200,11.250)(3.800,11.250)
\psline[linewidth=0.5mm]{->}(0.200,9.750)(3.800,0.750)
\psline[linewidth=0.5mm]{->}(0.200,8.250)(3.800,8.250)
\psline[linewidth=0.5mm]{->}(0.200,6.750)(3.800,5.250)
\psline[linewidth=0.5mm]{->}(0.200,5.250)(3.800,2.250)
\psline[linewidth=0.5mm]{->}(0.200,3.750)(3.800,9.750)
\psline[linewidth=0.5mm]{->}(0.200,2.250)(3.800,6.750)
\psline[linewidth=0.5mm]{->}(0.200,0.750)(3.800,3.750)
\end{pspicture}
\end{center}
\caption[xyz]{Depiction of the permutation
$\pi(1), \pi(2),\ldots, \pi(8)
= 1, 8, 3, 5, 7, 2, 4, 6$, which has $2$ descents.
Arrows that originate at $i$ and $i+1$ on the left  intersect
if $i=2$ or $i=5$.
If this permutation is to be
realized as an $a$-shuffle, cuts must be placed in the two positions shown
in Figure 1.  }
\label{fig-2}
\end{figure}

Every $a$-shuffle is a rearrangement of $n$ cards and can therefore be thought of as
a permutation $\pi$ of $\{1,2,\ldots,n\}$. If the $i$th card ends up in the $j$th position,
then $\pi(i) = j$. The number of {\it descents} of a permutation turns out to be an essential
concept and we will explain it in three different ways.

\begin{itemize}
\item Firstly, the number of descents
of $\pi$ is equal to the number of positions $i$, $1\leq i\leq n-1$,
with $\pi(i)>\pi(i+1)$. If we walk along the sequence $\pi(1),
\pi(2),\ldots,\pi(n)$ from beginning to end, a descent must be
recorded every time there is a decrease in value.
\item For the second explanation, the permutation $\pi$ must be represented as shown
in Figure \ref{fig-2}. The numbers $1$ to $n$ are listed twice with
the listing on the left corresponding to positions of the unshuffled
deck, or the {\it source deck}, and that on the right to positions of
the shuffled deck, or the {\it target deck}. If $\pi(i)=j$, an arrow
is drawn that originates at the $i$ on the left and terminates at the
$j$ on the right. The number of descents is equal to the number of
pairs of arrows that originate at consecutive positions $i$ and $i+1$
and cross each other.

\item For the third
explanation, we realize the permutation $\pi$ as an $a$-shuffle for
some $a$.  If $\pi$ is depicted as shown in Figure \ref{fig-2}, each
packet is a block of contiguous cards on the left, and the arrows
coming out of the same packet or block may not cross each other.
A cut must be placed wherever two arrows that originate at
consecutive positions cross each other. Therefore the minimum
number of cuts necessary to realize $\pi$ as an $a$-shuffle is
equal to the number of descents in $\pi$, and $a$ equals the
number of cuts plus $1$.
\end{itemize}

Bayer and Diaconis \cite{BD1} proved that
the probability that an $a$-shuffle results
in a permutation $\pi$, with $\des(\pi)=d$, is given by
\begin{equation}
\frac{1}{a^n}\binom{a+n-d-1}{n}.
\label{eqn-2-1}
\end{equation}
What we need to understand is the effect of repeated $2$-shuffles of a
deck. The above formula alone is not enough.  The result stating that
an $a$-shuffle followed by a $b$-shuffle is equivalent to an
$ab$-shuffle, proved by Aldous \cite{Aldous1}, Bayer and Diaconis
\cite{BD1}, and Reeds, is also needed.  With that result, it follows
immediately that $k$ $2$-shuffles are equivalent to a single
$2^k$-shuffle. The probability that $k$ $2$-shuffles result in a
permutation $\pi$ can then be found using \eqref{eqn-2-1}. In
addition, \eqref{eqn-2-1} tells us that the transition probability can
take on only $n$ different values corresponding to $d=0,1,\ldots, n-1$
for $a$ fixed.

\section{Decks with repeated cards and the descent polynomial}
In this paper, the term {\it deck} refers to an ordered sequence of cards.
Let $D$ be a deck with cards labeled $1, 2, \ldots, h$. If the number
of cards labeled $1$ is $n_1$, the number labeled $2$ is $n_2$, and so
on, the total number of cards in $D$ is $n=n_1+n_2+\cdots+n_h$.
For example, the
deck $D=1,2,1,2,1,1$ has $n_1=4$, $n_2=2$, and $n=6$.  The deck
\[D=\underbrace{1\ldots1}_{n_1}\underbrace{2\ldots 2}_{n_2},\ldots, 
\underbrace{h\ldots h}_{n_h},\]
which we will abbreviate as $1^{n_1},2^{n_2},\ldots, h^{n_h}$, has
$n_1$ cards labeled $1$ above $n_2$ cards labeled $2$ and so on. The
deck $D=(1,2)^{n_0}$ has $n=2n_0$ cards with cards labeled $1$ and $2$
alternating. Given $D$, $D(i)$ denotes the label of the card in the
$i$th position in $D$. For example, if $D = 1,2,1,2,1,1$ then $D(1) =
1$ and $D(4)=2$.

Normally, the cards of a deck will be listed from left to right, with
the label of the topmost card appearing first in the list, and with the
labels separated by commas, as in  the previous paragraph. However, sometimes
the commas will be omitted. 

If a permutation $\pi$ is applied to a deck $D=e_1,\ldots, e_n$, it
sends the card $e_i$ in position $i$ to position $\pi(i)$. Therefore
the resulting deck is $e_{\pi^{-1}(1)}, e_{\pi^{-1}(2)},\ldots,
e_{\pi^{-1}(n)}$, where $\pi^{-1}$ is the inverse of the permutation
$\pi$.

Let $D_1$ and $D_2$ be decks of $n$ cards, $n_c$ of which are labeled
$c$ for $1\leq c\leq h$. We say that a permutation $\pi$ of
$\{1,2,\ldots,n\}$ belongs to the set of permutations from $D_1$ to
$D_2$ if $D_1(i) = D_2(\pi(i))$ for $1\leq i \leq n$. That set will
be denoted by $\Pi(D_1; D_2)$. It
includes all the permutations which when applied to $D_1$ result in
$D_2$ and only those. For example, if $D_1 =
1,1,2,2$ to $D_2 = 1,2,2,1$, $\Pi(D_1; D_2)$
has $4$ members, given by $\pi(1), \pi(2),
\pi(3),
\pi(4)$ equal to $1,4,2,3$ or $1,4,3,2$ or $4,1,2,3$ or $4,1,3,2$. In
general the cardinality of $\Pi(D_1; D_2)$ is
$n_1!n_2!\ldots n_h!$.

The descent polynomial of $\Pi(D_1;D_2)$,
the set of permutations from $D_1$ to $D_2$, is defined as 
\( \sum_{\pi} x^{\des(\pi)},\)
where $\pi$ ranges over the set $\Pi(D_1;D_2)$ and
$\des(\pi)$ is the number of descents of $\pi$. Let the descent
polynomial be
\begin{equation}
c_0 + c_1 x +\cdots+c_{n-1} x^{n-1}.
\label{eqn-3-1}
\end{equation}
The coefficient $c_d$ equals the number of permutations in $\Pi(D_1;D_2)$
with $d$ descents.  The probability that an $a$-shuffle of $D_1$
results in the deck $D_2$ is therefore given by
\begin{equation}
p_a = \sum_{d=0}^{n-1}\frac{c_d}{a^n} \binom{a+n-d-1}{n},
\label{eqn-3-2}
\end{equation} 
a formula obtained by using \eqref{eqn-2-1} and summing over all
the permutations in $\Pi(D_1;D_2)$. Setting $a=1,2,\ldots,n$
gives a triangular system of equations for $p_1, p_2,\ldots, p_n$ in
terms of the coefficients $c_0, c_1,\ldots, c_{n-1}$.  This triangular
system can be inverted using a binomial identity (see
\cite[p. 269]{GKPBook}) to get
\begin{equation}
c_d = p_{d+1}(d+1)^n - p_d d^n \binom{n+1}{1}+
p_{d-1}(d-1)^n\binom{n+1}{2}-\cdots+(-1)^d p_1 1^n \binom{n+1}{d},
\label{eqn-3-3}
\end{equation}
for $1\leq d < n$. Using \eqref{eqn-3-2} and \eqref{eqn-3-3}, it is easy to pass
back and forth between the transition probabilities $p_a$ and the descent polynomial
\eqref{eqn-3-1}.

Suppose we are given a source deck $D_1$ and asked to determine how
many riffle shuffles mix that deck. We will examine a
definition of mixing later on, but surely we need to be able to
determine the transition probability from $D_1$ to any rearrangement
of its cards under an $a$-shuffle. Since it is easy to pass back and
forth between the transition probabilities and the descent polynomial,
we need to be able to determine the descent polynomial of the
permutations
from $D_1$ to any rearrangement of it. This begs the question, given decks
$D_1$ and $D_2$ is it possible to compute the descent polynomial
\eqref{eqn-3-1} of permutations from $D_1$ to $D_2$ efficiently? To answer this question, we take a trip
through computational complexity theory. 

\section{An excursion to computational complexity theory}

We begin with a decision problem seemingly unrelated to our concerns:
\vspace*{.2cm}

\noindent{THREE DIMENSIONAL MATCHING (3DM):}
Finite sets $X=\{x_1,\ldots, x_m\}$, $Y=\{y_1,\ldots,y_m\}$, and
$Z=\{z_1,\ldots,z_m\}$ of equal cardinality are given.  A subset $T$
of $X\times Y\times Z$ is also given. Decide if there exists a subset
$M$ of $T$ of cardinality $m$ such that every element of $X$ occurs
exactly once as the first element of a triple in $M$, every element of
$Y$ occurs exactly once as the second element of a triple in $M$, and
every element of $Z$ occurs exactly once as the third element of a
triple in $M$.

\vspace*{.2cm}
3DM belongs to the class of decision problems known as NP. Karp
\cite{Karp1} showed a way to reduce every problem in NP to 3DM in time
polynomial in problem size and thus proved 3DM to be NP-complete. Karp
proved the NP-completeness of a number of other decision problems as
well; his paper is a cornerstone of computational complexity theory.
The $\mathrm{P}\neq\mathrm{NP}$ conjecture implies that there is no
algorithm for 3DM whose running time is polynomial in $m$. Finding a
polynomial time algorithm or a super-polynomial lower bound for a
single NP-complete problem suffices to resolve this conjecture since
all NP-complete problems can be reduced to each other in polynomial
time. Kozen's book \cite{KozenBook} has a good introduction to some of
the basic topics of computational complexity such as NP-completeness
and \#P-completeness.

The following decision problem is related to card shuffling:

\vspace*{.2cm}
\noindent{MIN CUTS:} Given two decks, $D_1$ and $D_2$, and a positive integer $d$,
decide if there exists a permutation $\pi\in \Pi(D_1;D_2)$ 
with $d$ or fewer descents.

\vspace*{.2cm}
To determine the transition probability from $D_1$ to $D_2$ under an
$a$-shuffle using \eqref{eqn-3-2}, we need to know the coefficients
$c_i$ of the descent polynomial \eqref{eqn-3-1}. MIN CUTS asks if one
of the coefficients $c_i$, $0\leq i\leq d$, is nonzero for a given
$d$.

Another decision problem related to card shuffling is the following:

\vspace*{.2cm}
\noindent{RIFFLE:} Given $p$ nonempty packets of cards $P_1, P_2,\ldots, P_p$ and a deck $D$,
decide if it is possible to riffle the packets and get $D$. 

\vspace*{.2cm}
Each packet of cards $P_i$ is a set of cards with labels placed one
above the other.  We may also call it a deck, although we do not.  All
the cards in all the packets will appear somewhere in the new
deck. RIFFLE can also be worded in terms of sequences and
subsequences.  RIFFLE has a solution if and only if each packet $P_i$
corresponds to a subsequence of $D$, such that each packet equals the
corresponding subsequence and every card of $D$ appears in exactly one
subsequence.

Both RIFFLE and MIN CUTS are NP-complete. To prove that, we give
polynomial time many-one reductions from 3DM to RIFFLE and from RIFFLE
to MIN CUTS. 

\vspace*{.2cm}
\noindent{\bf 3DM reduces to RIFFLE.}
 Suppose we are given an instance of 3DM. It can be
reduced to RIFFLE as follows.  The elements $x_i, y_i, z_i$ of $X,Y,Z$
become card labels and we assume $X,Y,Z$ to be pairwise disjoint with
no loss of generality.  For every triple $(x_i, y_j, z_k) \in T$
create the packet $x_i, y_j, z_k, L$, where $L$ is a new card
label. For example, if $(x3, y1, z2)\in T$, we create the packet $x3,
y1, z2, L$. If there are $t$ triples in $T$, we create $t$ packets
totally in the instance of RIFFLE. The deck $D$ for this instance of
RIFFLE is given by
\[ D = x_1,\ldots,x_m,y_1,\ldots, y_m, z_1,\ldots, z_m, L^m, 
x_1^{\alpha_1},\ldots, x_m^{\alpha_m}, y_1^{\beta_1},\ldots, y_m^{\beta_m},
z_1^{\gamma_1},\ldots,z_m^{\gamma_m}, L^{m^\ast},\]
where $\alpha_i$, $\beta_i$, and $\gamma_i$ are one less than the number of
occurrences of $x_i$, $y_i$, and $z_i$ as elements of triples in $T$ , respectively,
but if the number of occurrences is zero those numbers must also be zero.
Additionally, $m^{\ast} = \max(0, t-m)$, where $t$ is the number of triples
in $T$. 
We claim that the given instance of 3DM has a matching $M$ if and only if
the $t$ packets in the created instance of RIFFLE can be riffled to get $D$.

Suppose the matching $M$ is a solution of 3DM and $(x_i, y_j, z_k)\in M$.
Then the card labeled $L$ in the packet $x_i, y_j, z_k, L$ must be dropped
to get one of the $m$ $L$s that occur as the first block of $L$s in $D$,
and the cards $x_i$, $y_j$, $z_k$ must be dropped from this packet to get
the cards with those labels that precede the first block of $L$s in $D$.
If the triple $(x_i,y_j,z_k)$ belongs to $T$ but not to $M$, the
card labeled $L$ in the packet $x_i,y_j,z_k,L$ must be dropped to get one of
the $m^{\ast}$ $L$s that occurs at the very bottom of $D$ and the cards 
$x_i$, $y_j$, $z_k$ must be dropped from this packet to get cards with those
labels that follow the first block of $L$s in $D$. We can riffle the packets 
in this way to get $D$. Suppose the instance of RIFFLE that was created from
the given instance of 3DM has a solution. Consider the $m$ packets whose
$L$s are dropped to get the $m$ $L$s that occur as the first block of $L$s
in $D$. The triples that correspond to these packets must form a matching
$M$ of 3DM.  The proof of validity of the reduction from 3DM to RIFFLE
is now complete.

The reduction from 3DM to RIFFLE uses cards with $3m+1$ different
labels. What if cards are not allowed to have arbitrarily many labels?
We modify the reduction to use only cards with the four labels $[$,
$]$, $c$, and $L$.  The modification is to replace each $x_i$ that
occurs in a packet or in the deck $D$ by the list of cards $[ c^i ]$;
each $y_i$ by $[ c^{m+i} ]$; and each $z_i$ by $[ c^{2m+i} ]$. The
opening and closing brackets $[$ and $]$ that occur in $D$ and in
each of the packets are properly matched with no nesting. Therefore,
in any solution of RIFFLE, if a packet drops a $[$ to get the very
last or bottommost $[$ that occurs in $D$, that dropped $[$ must be
the bottommost $[$ in that packet and the same packet must have
dropped the matching $]$ to get the bottommost $]$ in $D$.  By
induction, we conclude that matching parentheses $[$ and $]$ in a
packet must be dropped to get matching parentheses in $D$. Therefore
the list of cards that code for $x_i$ or $y_i$ or $z_i$ in any packet
must be dropped all at once to get a contiguous sequence of cards in
$D$.  The argument for the validity of the earlier reduction can now
be applied to its modification.  A little thought will convince the
reader that cards labeled $L$ in this reduction can be replaced by
cards labeled $c$.  We have proved that RIFFLE is NP-complete even if
the cards are allowed to have only three different labels.

\vspace*{.2cm}
\noindent{\bf RIFFLE reduces to MIN CUTS.}
RIFFLE can be reduced to MIN CUTS as follows. Given an instance 
$(P_1,\ldots, P_p;D)$ of RIFFLE,
consider the decks
\begin{align*}
D_1 &= P_1 L P_2 L \ldots P_p \\
D_2 &= D L^{p-1},
\end{align*}
where $L$ is a new label that does not occur in the given instance of RIFFLE.
The instance of MIN CUTS uses these decks $D_1$ and $D_2$ and chooses $d=p-1$.
We claim that RIFFLE has a solution if and only if there exists a permutation
in $\Pi(D_1;D_2)$ with $p-1$ descents. For a proof, it
suffices to note that any permutation in $\Pi(D_1;D_2)$ must
cut $D_1$ after each of the $p-1$ $L$s in $D_1$ as all the $L$s in $D_2$ occur at
the bottom. This reduction uses only a single new label. We have now proved the
following theorem.

\begin{thm}
MIN CUTS is NP-complete. MIN CUTS remains NP-complete even if cards are allowed
only four different labels.
\label{thm-4-1}
\end{thm}

The problem of finding the descent polynomial \eqref{eqn-3-1} of 
$\Pi(D_1; D_2)$
is the counting version of MIN CUTS since the coefficient $c_d$ equals the
number of permutations in $\Pi(D_1;D_2)$ with $d$ descents. Once the definition of \#P-completeness
in terms of Turing machines is understood, it will  be clear that the proof of 
Theorem \ref{thm-4-1}, with slight modifications, implies that
finding the descent polynomial is \#P-complete
even if only four different labels are allowed. We conjecture that
finding the descent polynomial is \#P-complete even if $D_1$ and $D_2$ 
have cards with only two different labels. 

If all cards in $D_1$ are distinct, there can be only one permutation $\pi$ in
$\Pi(D_1;D_2)$ and the descent polynomial can be easily found. There are other
interesting cases where the descent polynomial can be found efficiently.
If either $D_1$ or $D_2$ is a deck where all cards with the same label occur
in a sequence of contiguous positions, there is an efficient algorithm for
finding the descent polynomial. See \cite{CV1}. There may be yet more cases.

There is an algorithm of time complexity $O(n^{2d})$ for determining the
coefficient of $x^d$ in the descent polynomial of permutations from a given
deck $D_1$ to another given deck $D_2$. A detailed description of this algorithm
will appear in the first author's Ph.D. thesis. As the running time of this
algorithm is exponential in $d$, its use to determine the mixing times for
games like blackjack and bridge is impractical on today's computers.

\begin{figure}
\begin{center}
\includegraphics[height=2in,width=2in]{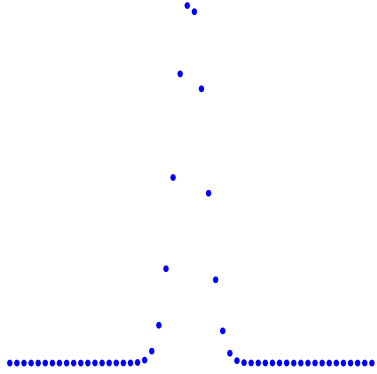}
\end{center}
\caption[xyz]{A graph of the coefficients of a descent polynomial
looks like the bell curve.}
\label{fig-X}
\end{figure}

\section{A brief introduction to Stein's method of auxiliary randomization}
As the task of finding the descent polynomial is \#P-complete and
therefore unlikely to have an efficient algorithm, it is necessary to
find a way to approximate the descent polynomial.

 To approximate the
descent polynomial, we turn to Stein's method of auxiliary
randomization. Stein's method can be used to approximate some
probability distributions by the normal distribution or the Poisson
distribution or the binomial distribution.

Graphs of coefficients of most descent polynomials look like bell curves,
as shown in Figure \ref{fig-X} and as proved in the next section.
Although the normal approximations give us a good sense for what
these polynomials look like, the approximations are not accurate for
the first few coefficients of the descent polynomial. We will need
those coefficients with good relative accuracy, and in Section 7,
we obtain usable approximations for those coefficients based on
heuristic arguments. The methods of Section 7 are guided by 
normal approximation results.

Let $X_1, X_2,\ldots, X_n$ be a sequence of independent and identically
distributed random variables with $\P(X_i=0) = \P(X_i=1) = 1/2$. 
By elementary probability theory, the distribution
of the sum $S=X_1+X_2+\cdots + X_n$ is binomial and can be approximated by the
standard normal distribution after suitable normalization. We use this simple
example to give a brief introduction to Stein's method \cite{Stein1} \cite{SteinBook}.
The introduction is broken down into three steps.

\vspace*{.2cm}
\noindent{\bf Characterization of the standard normal distribution.}
The standard normal distribution has the probability density function
$(1/\sqrt{2\pi}) \exp(-x^2/2)$. Another characterization would be in terms of
its moments. There are ways to prove central limit theorems using these characterizations
of the normal distribution, but neither is used by Stein's method. 
The characterization of the normal distribution used by Stein's method is given
in the lemma below.
 
\begin{lem}
A random variable $W$ has the standard normal distribution if and only if the
expectations $\E(Wf(W))$ and $\E(f'(W))$ are equal for all bounded continuous functions $f$
whose first and second derivatives are bounded and piecewise continuous.
\label{lem-4-1}
\end{lem} 

The proof in one direction is a simple use of integration by
parts. The proof in the other direction is also elementary, but more
involved. Stein \cite[chapter 2]{SteinBook} has given a proof of a much refined
version of this lemma.

For the sum $S$, $\mu = \E S = n/2$ and $\sigma^2 = \var(S) = n/4$. We want to show that the distribution
of $T = (S-\mu)/\sqrt{\sigma}$ is close to normal. The expectations $\E(Tf(T))$
and $\E(f'(T))$, for functions $f$ as in Lemma \ref{lem-4-1}, will not be exactly equal,
since the distribution of $T$ is not exactly normal, but if it is possible to show
that the expectations are close to one another then there might be a way to show that
the distribution of $T$ is close to normal. But how to show that the  expectations $\E(Tf(T))$
and $\E(f'(T))$ are close? The answer is to proceed by throwing some extra randomness
into the problem and then use the familiar Taylor series expansion.

\vspace*{.2cm}
\noindent{\bf Auxiliary randomization.}
There are many ways to add extra randomness with
a view to showing the expectations $\E(Tf(T))$ and $\E(f'(T))$ to be close to each
other. The way described here is related to size biasing. 

To explain size biasing, we consider a random variable $W$ with a
continuous density function $p(x)$ and finite expectation. Assume that
$W\geq 0$ or $p(x)=0$ for $x<0$. Then $\E W = \int_0^\infty x p(x)
dx$. Therefore $xp(x)/\E W$ is also a probability density function. A
variable with that density is said to be $W$-size biased. More
generally, $W^\ast$ is said to have the $W$-size biased distribution
if $\E(Wf(W)) = \E W \E(f(W^\ast))$ for all continuous $f$ for which
the expectation on the left hand side exists. The size biased
distribution is defined only for non-negative random variables. If $W$
is a $0$-$1$ valued random variable which is $1$ with a positive
probability, its size biased distribution assigns probability $1$ to
the value $1$.
 
The size biasing of sums such as $S$ above is facilitated by the following lemma
\cite{BRS1}. 
\begin{lem}
Let $W = X_1 + X_2+\cdots+X_n$ where the $X_i$ are all $0$-$1$ valued random variables,
but not necessarily independent or identically distributed.  
Let $I$ be a random variable which is independent of the $X_i$ and which satisfies
$\P(I=i) = \E X_i/\E W$. Let $X_1^\ast, X_2^\ast, \ldots, X_n^\ast$ be a sequence
of random variables such that $X_I^\ast = 1$ and
\[\P\bigl((X_1^\ast, X_2^\ast, \ldots, X_n^\ast)\in A\bigl | I=i\bigr)
= \P\bigl((X_1, X_2,\ldots, X_n)\in A\bigl | X_i = 1\bigr),\]
for all possible sets $A$. Then $W^\ast = X_1^\ast+X_2^\ast+ \cdots+ X_n^\ast$ has
the $W$-size biased distribution
\label{lem-5-2}
\end{lem}
\begin{proof}
The sequence of equalities below proves the lemma.
To justify the second equality below, note the assumption about the conditional
distributions of the $X_i^\ast$ and the $X_i$ in the lemma. For the fourth equality
below, note that $\E X_i =\P (X_i = 1)$ since $X_i$ is $0$-$1$ valued.
\begin{align*}
\E f(W^\ast) &= \sum_{i=1}^n \E\bigl(f(W^\ast)\bigl | I=i\bigr) \P\bigl(I=i\bigr) \\
&= \sum_{i=1}^n
\E\bigl(f(W)\bigl | X_i=1\bigr) \P\bigl(I=i\bigr)\\
&= \sum_{i=1}^n \E\bigl(f(W)\bigl | X_i=1\bigr) \E X_i/ \E W\\
&= \frac{1}{\E W}\sum_{i=1}^n \E\bigl(f(W)\bigl | X_i=1\bigr) \P\bigl( X_i=1\bigr)\\
&= \frac{1}{\E W} \sum_{i=1}^n \E\bigl(X_i f(W)\bigr) = \frac{\E(W f(W))}{\E W}.
\end{align*}
\end{proof}

Recall that the random variable $S$ was defined to be a sum of
independent and identically distributed random variables
$X_1,\ldots,X_n$. Let $I$ be independent of the $X_i$ and uniformly
distributed over the set $\{1,2,\ldots, n\}$. Define $X_i^\ast = X_i$
if $i\neq I$ and $X_I^\ast = 1$. Now the hypotheses of Lemma
\ref{lem-5-2} are easily verified and we may assert that 
the random variable $S^\ast =
X_1^\ast+\cdots+X_n^\ast$ has the $S$-size biased distribution.

The random variable $S^\ast$ was built  up using the $X_i$ and some extra randomness,
namely the random variable $I$. Taylor series expansion must be used to see how
having $S^\ast$ around helps in proving the distribution of $S$ to be normal.

\vspace*{.2cm}
\noindent{\bf Taylor series expansion.} 
Let $W$ be a non-negative random variable. Let $W^\ast$ have the $W$-size biased
distribution. Assume $\mu = \E W$ and $\sigma^2 = \var(W) = \E(W-\mu)^2$.  We want
to compare the distribution of $(W-\mu)/\sigma$ with the standard normal distribution. Consider
the following calculation.
\begin{align*}
\E\Biggl(\biggl(\frac{W-\mu}{\sigma}\biggr) f\biggl(\frac{W-\mu}{\sigma}\biggr)\Biggr)
&= \E\Biggl(\frac{W}{\sigma}f\biggl(\frac{W-\mu}{\sigma}\biggr)\Biggr)
- \frac{\mu}{\sigma} \E f\biggl(\frac{W-\mu}{\sigma}\biggr)\\
&= \frac{\mu}{\sigma}\Biggl( \E f\biggl(\frac{W^\ast-\mu}{\sigma}\biggr)
- \E f\biggl(\frac{W-\mu}{\sigma}\biggr)\Biggr)\\
&= \frac{\mu}{\sigma} \E\Biggl( \frac{(W^{\ast}-W)}{\sigma} f'\biggl(\frac{W-\mu}{\sigma}
\biggr)\Biggr) \\ &+
\frac{\mu}{\sigma} \E\Biggl( \frac{(W^\ast-W)^2}{\sigma^2}
f''\biggl(\frac{W-\mu}{\sigma}+\biggl(\frac{W^\ast-W}{\sigma}\biggr)U\biggr)(1-U)\Biggr),
\end{align*}
where $U$ is uniformly distributed over $[0,1]$ and independent of all
other random variables.  In the last line above, the remainder term
from Taylor's theorem, which is usually written as an integral, is
written as an expectation by using $U$. A simple calculation shows
that $\E (W^\ast - W) = (\E W^2/\E W) - EW = \sigma^2/\mu$. This
calculation makes it plausible that the first term on the right hand
side of the last equality above could be close to $\E
f'((W-\mu)/\sigma)$. The second term will be small only if $W^\ast-W$
is small. This is an important requirement. Not only must $W^\ast$
have the $W$-size biased distribution, the joint distribution of $W$
and $W^\ast$ must be such that $W^\ast$ is close to $W$.

Estimates carried out after Taylor series expansion imply the following theorem. It
is due to Goldstein \cite{Goldstein1}. The size biased version of Stein's method was
derived by Baldi, Rinott and Stein \cite{BRS1}.


\begin{thm}
Let $W$ be a non-negative random variable with $\E W = \mu$ and $\var(W) = \sigma^2$.
Let $W^\ast$ be jointly defined with $W$ such that its distribution is $W$-size
biased. Let $\abs{W^\ast-W} \leq B$ and let $A = B/\sigma$. 
Let $B \leq \sigma^{3/2}/\sqrt{6\mu}$.
Then
\[\ABS{\P\Biggl(\frac{W-\mu}{\sigma} \leq x\Biggr) - \Phi(x)}
\leq 0.4 A + \frac{\mu}{\sigma}(64A^2+4A^3)+\frac{23\mu}{\sigma^2}
\sqrt{\var\bigl(\E(W^\ast-W|W)\bigr)},\]
where $\Phi$ is the standard normal distribution.
\label{thm-5-3}
\end{thm}

Using Theorem \ref{thm-5-3}, we can prove that the distribution of $(S-\mu)/\sigma$
is close to the normal distribution. We state again that $\mu = \E S = n/2$
and $\sigma^2 = \var(S) = n/4$. From the construction of $S^\ast$, 
$\abs{S^\ast- S} \leq B$ for $B=1$. Given $S$, $S^\ast-S=1$ if $I=i$ and 
$X_i=0$. This happens with probability $(n-S)/n$. Otherwise, $S^\ast - S = 0$.
Therefore, $\E\bigl(S^\ast-S\bigl | S\bigr) = 1 -\frac{S}{n}$ and 
$\var\bigl(\E\bigl(S^\ast-S\bigl | S\bigr)\bigr) = 1/(4n)$. By Theorem \ref{thm-5-3},
\[\ABS{\P\Biggl(\frac{S-\mu}{\sigma} \leq x\Biggr) - \Phi(x)} \leq \frac{C}{\sqrt{n}},\]
for some constant $C$.

\section{Normal approximations for descent polynomials}
A number of normal approximation theorems for descents and inversions of permutations
of sets and multisets can be proved using the size biased version of Stein's method
\cite{CV2}. In this section, we obtain a normal approximation for the coefficients of
the descent polynomial \eqref{eqn-3-1} in a quite general setting. We use the technique
developed in \cite{CV2}.

Let $\pi$ be a uniformly distributed permutation in $\Pi(D_1;D_2)$,
where $D_2$ is obtained by rearranging the cards of $D_1$ in some
order.  Let the decks have $n_c$ cards with label $c$ for $1\leq c
\leq h$ and let the total number of cards be $n$.  The random
variables $X_1, X_2, \ldots, X_{n-1}$ are defined as follows: $X_i =
1$ if $\pi(i)>\pi(i+1)$ but $X_i=0$ otherwise. Let $W = \des(\pi) = 
X_1+X_2+\cdots+X_{n-1}$. Then
\begin{equation} \P(W=d) = \frac{c_d}{n_1!n_2!\ldots n_h!}
\label{eqn-6-1},
\end{equation}
where $c_d$ is the coefficient of the $x^d$ term of the descent polynomial \eqref{eqn-3-1}.
If an approximation for the distribution of $W$ is available, \eqref{eqn-6-1} can be   
used to approximate $c_d$.

We will construct $W^\ast$ so that its distribution is $W$-size biased.
Other constructions of this type can be found in \cite{CV2}.
 Let $I$ be a random
variable independent of $\pi$ with $\P(I=i) = \E X_i/\E W$ for 
$1\leq   i\leq n-1$.
We assume $\E W > 0$ so that $I$ is well defined.
If $c$ and $e$
are two distinct card labels it is useful to define the following
set:
\begin{equation}
S(c,e) = \bigl\{ (x,y) \bigl| x>y, D_2(x)=c, D_2(y) = e\bigr\},
\label{eqn-6-2}
\end{equation}
where $D(i)$, as noted earlier, stands for the label of the $i$th card in the deck $D$. 
Note $\E X_i = \Abs{S(c,e)}/(n_c n_e)$ if $D_1(i) = c$, $D_1(i+1) = e$ and $c\neq e$,
but $\E X_i = 1/2$ if $D_1(i) = D_1(i+1)$.

Given $\pi$ and $I$, we define a permutation
$\pi^\ast$, with $\pi^\ast\in\Pi(D_1;D_2)$. If $\pi(I) > \pi(I+1)$, then $\pi^\ast = \pi$.
If $\pi(I) < \pi(I+1)$ and $D_1(I) = D_1(I+1)$, then $\pi^\ast(I) = \pi(I+1)$,
$\pi^\ast(I+1) = \pi(I)$, and $\pi^\ast(i) = \pi(i)$ for $i\neq I, I+1$. The
remaining case is $\pi(I) < \pi(I+1)$, $D_1(I) = c$, $D_1(I+1)=e$, and $c\neq e$.
Let $J$ be a random variable independent of $\pi$ and $I$, and 
uniformly distributed over the set $S(c,e)$ of \eqref{eqn-6-2}.
Let $J=(x,y)$. 
Consider the list
$\pi(1), \pi(2),\ldots, \pi(n-1)$. 
Exchange $\pi(I)$ and $x$ and exchange
$\pi(I+1)$ and $y$ to get a new list. The permutation $\pi^\ast$ from $D_1$ to $D_2$ is
defined by setting $\pi^\ast(1),\pi^\ast(2),\ldots,\pi^\ast(n-1)$ equal to this
new list.

There is another way to describe $\pi^\ast$ in the last case above, which is
$\pi(I) < \pi(I+1)$, $D_1(I) = c$, $D_1(I+1)=e$, and $c\neq e$. Define
\begin{equation}
S_\pi(c,e) = \bigl\{(k,l)\bigl| (\pi(k),\pi(l))\in S(c,e)\bigr\}.
\label{eqn-6-3}
\end{equation}
We can pick $(k,l)$ uniformly from this set, exchange $\pi(I)$ with $\pi(k)$,
and exchange $\pi(I+1)$ with $\pi(l)$ to get $\pi^\ast$. The set over which
$(k,l)$ is distributed depends upon $\pi$ and $I$, but the distribution 
is always uniform.

\begin{lem}
$\P\bigl(\pi^\ast \in A\bigr| I = i \bigr) = \P\bigl(\pi \in A\bigl| \pi(i) > \pi(i+1)\bigr)$
for any possible set $A$.
\label{lem-6-1}
\end{lem}
\begin{proof}
If $I=i$ and $\pi(i) > \pi(i+1)$, then $\pi^\ast = \pi$. Therefore it suffices to show
that
\[\P\bigl(\pi^\ast \in A\bigr| I = i, \pi(i) < \pi(i+1) \bigr) = \P\bigl(\pi \in A\bigl| \pi(i) > \pi(i+1)\bigr).\]
The first case is when $D_1(i) = D_1(i+1)$. In this case, we have
\begin{align*}
&\P\bigl(\pi^\ast \in A\bigr| I = i, \pi(i) < \pi(i+1) \bigr)\\ =
&\P\bigl(\pi(1),\ldots,\pi(i-1),\pi(i+1),\pi(i),\pi(i+2),\ldots,\pi(n-1)\in A\bigl|
I = i, \pi(i) < \pi(i+1)\bigr)\\
=&\P\bigl(\pi(1),\ldots,\pi(i-1),\pi(i+1),\pi(i),\pi(i+2),\ldots,\pi(n-1)\in A\bigl|
\pi(i) < \pi(i+1)\bigr)\\
=& \P\bigl(\pi \in A\bigl| \pi(i) > \pi(i+1)\bigr)
\end{align*}
The last equality above holds because if we pick a uniformly distributed permutation
in $\Pi(D_1;D_2)$ and exchange $\pi(i)$ and $\pi(i+1)$, the new permutation is
also a uniformly distributed permutation in $\Pi(D_1;D_2)$ provided $D_1(i)
= D_1(i+1)$.

The second case is when $D_1(i) = c$, $D_1(i+1) = e$, and $c\neq e$. 
Below each of the summations is taken over $(x,y)\in S(c,d)$.
\begin{align*}
&\P\bigl(\pi^\ast\in A \bigl| I = i, \pi(i) < \pi(i+1)\bigr)\\
= &\sum \P\bigl(\pi^\ast\in A\bigl| I = i, \pi(i) < \pi(i+1), J = (x,y)\bigr)
\P\bigl(J=(x,y)\bigr)\\
= &\sum \P\bigl(\pi^\ast\in A\bigl| I = i, \pi(i) < \pi(i+1), J = (x,y)\bigr)
\P\bigl(\pi(i)=x, \pi(i+1)=y\bigl| \pi(i)>\pi(i+1)\bigr)\\
= &\sum\P(\pi \in A\bigl| \pi(i) = x, \pi(i+1)=y\bigr)
\P\bigl(\pi(i)=x, \pi(i+1)=y\bigl| \pi(i)>\pi(i+1)\bigr)\\
= &\P\bigl(\pi \in A\bigl| \pi(i) > \pi(i+1)\bigr).
\end{align*}
The second equality above holds because $J$ is uniformly distributed over the
set $S(c,d)$ of \eqref{eqn-6-2} and because $\pi$ is a uniformly distributed
permutation in $\Pi(D_1;D_2)$. The third equality above holds because $\pi$
is a uniformly distributed permutation in $\Pi(D_1;D_2)$ and because of the
way $\pi^\ast$ is generated  using $\pi$, $I$, and $J$.
\end{proof} 

For $i=1,2,\ldots,n-1$, define $X_i^\ast = 1$ if $\pi^\ast(i) > \pi^\ast(i+1)$,
but $X_i^\ast = 0$ otherwise. Let $W^\ast = \des(\pi^\ast) = X_1^\ast + X_2^\ast+\cdots+X_{n-1}^\ast$.

\begin{lem}
$W^\ast$ has the $W$-size biased distribution.
\label{lem-6-2}
\end{lem}
\begin{proof}
Follows from Lemmas \ref{lem-5-2} and \ref{lem-6-1}.
\end{proof}

If Theorem \ref{thm-5-3} is to be applied to approximate the
distribution of $W$, it is necessary to find an upper bound for
$\var\bigl( \E(W^\ast-W|W)\bigr)$.  By \cite[p. 477]{BillingsleyBook},
$\var\bigl( \E(W^\ast-W|W)\bigr) \leq
\var\bigl( \E(W^\ast-W|\pi)\bigr)$, since $W$ is a function of $\pi$. Let
\begin{equation}
Q = \E(W^\ast-W|\pi).
\label{eqn-6-4}
\end{equation}
We will upper bound $\var(Q)$.

The identity
\begin{equation}
\var \Bigl(\sum_{i=1}^n Y_i\Bigr) =  \Biggl(\sum_{i=1}^n \var(Y_i)
+ 2 \sum_{1\leq i<j\leq n}\covar(Y_i,Y_j)\Biggr)
\label{eqn-6-5}
\end{equation}
holds because $\var(Y) = \E Y^2 - \bigl(\E Y\bigr)^2$ and
$\covar(Y, Z) = \E YZ - \E Y \E Z$. Along with \eqref{eqn-6-5}, the following
Lemma \ref{lem-6-3} is
useful for upper bounding $\var(Q)$.

The argument to upper bound $\var(Q)$ simplifies a great deal if it is assumed that
cards with any label occur the same number of times in $D_1$ or $D_2$; in other words,
$n_1=n_2=\cdots=n_h=n_0$ with $n_0\geq 1$. We make this assumption from here onwards
and put it to use immediately in the lemma below. With this assumption, the total
number of cards is $n=hn_0$.

In the lemma below, as will become evident from its proof,
the constants $7$ and $10$ can be replaced by
other positive integers. We take the constants as $7$ and $10$ to
facilitate later use of the lemma. The constants were chosen
to simplify the exposition and are not the best possible.

\begin{lem}
Let $\pi$ be a uniformly distributed permutation in $\Pi(D_1;D_2)$. Let
$f(\pi)$ depend only upon the relative order of 
$\pi(i_1), \pi(i_2),\ldots, \pi(i_r)$. Let $g(\pi)$ depend only upon the
relative order of 
$\pi(j_1), \pi(j_2),\ldots, \pi(j_s)$.
Assume $\abs{f}\leq 7$, $\abs{g}\leq 7$, $r\leq 10$, and $s\leq 10$.
There are three cases:
\begin{enumerate}
\item $\{i_1, \ldots, i_r\} \cap \{j_1,\ldots,j_s\} \neq \phi$,
\item $\{i_1, \ldots, i_r\} \cap \{j_1,\ldots,j_s\} = \phi$, but
$D_1(i_\rho) = D_1(j_\psi)$ for some $1\leq \rho \leq r$ and $1\leq \psi\leq s$,
\item $\{i_1, \ldots, i_r\} \cap \{j_1,\ldots,j_s\} = \phi$, and
$D_1(i_\rho) \neq D_1(j_\psi)$ for any $1\leq \rho \leq r$ and $1\leq \psi\leq s$.
\end{enumerate}
In these three cases, we have
\begin{enumerate}
\item $\covar(f,g)\leq C$ for some constant $C$ that depends on neither
$n_0$ nor $h$,
\item $\covar(f,g)\leq C/n_0$ for some constant $C$ that depends 
on neither $n_0$ nor $h$,
\item $\covar(f,g)=0$,
\end{enumerate}
respectively.
\label{lem-6-3}
\end{lem}
\begin{proof}
For the first case, it is enough to note that $f$ and $g$ are bounded in magnitude.
For the third case, it is enough to note that $f$ and $g$ are independent of each
other. The second case remains to be proved.

It is enough to consider $f$ and $g$ to be indicator functions that are $1$ for 
a particular relative ordering and $0$ for others. This is because any  $f$
or $g$ that is bounded in magnitude by a constant 
can be written as a linear combination of a constant number of such indicator
functions with coefficients that are bounded in magnitude by a constant.
We show the proof assuming $f=1$ if $\pi(i) < \pi(i+1)$ and $f=0$ otherwise;
and $g=1$ if $\pi(j) < \pi(j+1)$ and $g=0$ otherwise. For the second case
to apply, either $i+1< j$ or $j+1<i$ must hold, and 
at least one of $D_1(i), D_1(i+1)$ must equal
one of $D_1(j), D_1(j+1)$.

Let $P_1 = \P\bigl(\pi(i) < \pi(i+1)\bigr)$ and $P_2 = \P\bigl(\pi(j) < \pi(j+1)\bigr)$.
We claim that 
\begin{equation}
\P\bigl(\pi(j) < \pi(j+1)\bigl|\pi(i) = x, \pi(i+1)=y\bigr) = P_2 + \epsilon,
\label{eqn-6-6}
\end{equation}
with $\abs{\epsilon} < 4/n_0$. In \eqref{eqn-6-6}, $1\leq x,y\leq n$, $x\neq y$,
$D_2(x) = D_1(i)$, and $D_2(y) = D_1(i+1)$. The claim is true because of the following
argument. If $D_1(j) = D_1(j+1)$, the pair $(\pi(j), \pi(j+1))$ can take $n_0(n_0-1)$
different values. Otherwise, it can take $n_0^2$ values. For some of these values,
$\pi(j) < \pi(j+1)$. Given $\pi(i)=x$ and $\pi(i+1)=y$, at least one of $x$ or $y$
is not allowed to appear in a possible value for $(\pi(j), \pi(j+1))$. Thus at most
$2n_0$ possible values for this pair must be excluded.

The proof for this $f$ and $g$ can be completed by noting that $\covar(f,g)$ equals
\[\P\bigl(\pi(i)<\pi(i+1), \pi(j)<\pi(j+1)\bigl)
-\P\bigl(\pi(i)<\pi(i+1)\bigr)\P\bigl(\pi(j)<\pi(j+1)\bigr),\] 
and by writing the first of the three probabilities in the line above
in terms of the conditional probabilities on the left hand side of
\eqref{eqn-6-6}. The proof for
more general $f$ and $g$ is similar.
\end{proof} 

We now go back to $Q$ defined by \eqref{eqn-6-4} and the construction
of the size biased random variable $W^\ast$. Let us suppose that $\pi$ is a given
permutation in $\Pi(D_1;D_2)$. To write $Q$ as a sum, we introduce the quantities
$\chi_\pi(i,i+1)$ and $\chi_\pi(i,i+1, k, l)$, where $k\neq l$. The first of these 
is defined as the change in the number of descents when $\pi(i)$ and $\pi(i+1)$
are exchanged. The second is defined as the change in the number of descents
when $\pi(i)$ is exchanged with $\pi(k)$ and $\pi(i+1)$ is exchanged with $\pi(l)$.
With a view to applying Lemma \ref{lem-6-3} later on, we note that $\chi_\pi(i,i+1)$
depends only on the relative order of at most $4$ numbers, namely
$\pi(i-1), \pi(i), \pi(i+1), \pi(i+2)$. Similarly, $\chi_\pi(i,i+1,k,l)$ depends only on
the relative order of at most $10$ numbers of the form $\pi(i)$. The magnitude of
both of these quantities is always bounded by $7$.  

From here onwards, we denote $\E W$ and $\var(W)$ by $\mu$ and $\sigma^2$, respectively.

Let $A$ be the set of values of $i$ for which $\pi(i)<\pi(i+1)$
and $D_1(i) = D_1(i+1)$. Let $B$ be the set of values of $i$ for
which $\pi(i) < \pi(i+1)$,  $D_1(i)\neq D_1(i+1)$, and
$S(D_1(i),D_1(i+1))\neq \phi$. We must have
$A\cap B = \phi$. 
By the construction of $W^\ast$ using $\pi$, $I$ and $J$, if $i\in A$, we have
\begin{equation}
\P(I=i) = \frac{1}{2\mu} \quad \text{and}\quad 
\E\bigl(W^\ast-W\bigl|\pi, I=i\bigr) = \chi_\pi(i,i+1).
\label{eqn-6-7}
\end{equation}  
If $i\in B$, let $D_1(i) = c$ and $D_1(i+1) = e$. We have,
\begin{equation}
\P(I=i) = \frac{\Abs{S(c,e)}}{n_0^2}\frac{1}{\mu}
\label{eqn-6-8}
\end{equation}
and
\begin{align}
\E\bigl(W^\ast-W\bigl|\pi, I=i\bigr) &= \sum\E\bigl(W^\ast-W\bigl|\pi, I=i, J=(x,y)\bigr)
\P(J=(x,y))\nonumber\\
&= \sum\frac{\chi_\pi(i,i+1,k,l)}{\Abs{S(c,e)}},
\label{eqn-6-9}
\end{align}
where the first summation is over all $(x,y)$ in the set $S(c,e)$ of
\eqref{eqn-6-2} and the second
summation is over all $(k,l)$ in the set $S_\pi(c,e)$ of \eqref{eqn-6-3}. If
$i\notin A$ and $i\notin B$, either $P(I=i) = 0$ or $\E\bigl(W^\ast-W\bigl|\pi,I=i)=0$.

By \eqref{eqn-6-4}, \eqref{eqn-6-7}, \eqref{eqn-6-8} and \eqref{eqn-6-9}, we have
\begin{equation*}
Q = \sum_{i=1}^{n-1} \E\bigl(W^\ast-W\bigl|\pi, I=i\bigr) \P\bigl(I=i\bigr)
\end{equation*}
and
\begin{equation}
\mu Q = \sum_{i\in A} \frac{\chi_\pi(i,i+1)}{2}
+ \sum_{i\in B\;\text{and}\;(k,l)\in S_\pi(D_1(i),D_1(i+1))}
\frac{\chi_\pi(i, i+1, k, l)}{n_0^2}. 
\label{eqn-6-10}
\end{equation}
We will use \eqref{eqn-6-10} with \eqref{eqn-6-5} to upper bound $\var(\mu Q)$.

If \eqref{eqn-6-5} and \eqref{eqn-6-10} are used to write $\var(\mu Q)$ as a sum
of variance and covariance terms, each term on the right hand side of \eqref{eqn-6-10}
will contribute a single variance term. The total contribution of these variance
terms will be bounded by $Cn$, for some constant $C$, because each $\chi_\pi$ is
bounded in magnitude by $7$.

Some of the covariance terms will be of the form
\begin{equation}
\covar\Biggl(\frac{\chi_\pi(i_1,i_1+1,k_1,l_1)}{n_0^2},
\frac{\chi_\pi(i_2,i_2+1,k_2,l_2)}{n_0^2}\Biggr),
\label{eqn-6-11}
\end{equation}
where $i_1\in B$ and $i_2\in B$. The value of $\chi_\pi(i_1, i_1+1,k_1,l_1)$ depends
only upon the relative order of the $\pi(i)$, with $i$ equal to $i_1$ or $i_1+1$ or
$k_1$ or $l_1$, or differing from one of those $4$ integers by at most $1$.  There
can be at most $10$ such values for $i$ and we will denote this set of values
by $\nghd{i_1,i_1+1,k_1,l_1}$. As for the magnitude of the covariance term 
\eqref{eqn-6-11}, there are three cases.
\begin{enumerate}

\item Suppose $\nghd{i_1,i_1+1,k_1, l_1}\cap \nghd{i_2,i_2+1,k_2,l_2}\neq \phi$.
This corresponds to the first case of Lemma \ref{lem-6-3}. Therefore the magnitude
of the covariance term \eqref{eqn-6-11} is bounded by $C/n_0^4$ in this case.

To count the number of covariance terms of this type, note that $i_1$ can take at most
$n-1$ different values. As we require $D_1(k_1) = D_1(i_1)$ and $D_1(l_1)
= D_1(i_1+1)$, for given $i_1$, there are at most $n_0^2$ possible values of
$(k_1,l_1)$. For the covariance term \eqref{eqn-6-11} to fall under this type,
at least one of $i_2,k_2,l_2$ must differ from one of $i_1,k_1, l_1$ by less than
$2$. Therefore, given $i_1,k_1,l_1$, there are at most a constant number of choices
for one of $i_2,k_2,l_2$. Having chosen one of $i_2, k_2, l_2$, there are at most
$n_0^2$ ways to choose the other two. For example, suppose $k_2$ has been chosen.
This restricts $i_2$ to at most $n_0$ possibilities since we require
$D_1(i_2) = D_1(k_2)$. Given $i_2$, there are at most $n_0$ possibilities for
$l_2$ since we require $D_1(l_2) = D_1(i_2+1)$.  Thus the number of covariance terms
of this type is bounded by $C n n_0^4$ for some constant $C$.

\item Suppose $\nghd{i_1, i_1+1, k_1, l_1}\cap \nghd{i_2,i_2+1,k_2,l_2} = \phi$, but
$D_1(i^\ast_1) = D_1(i^\ast_2)$ for some $i^\ast_1 \in \nghd{i_1,i_1+1,k_1,l_1}$
and some $i^\ast_2\in\nghd{i_2,i_2+1,k_2,l_2}$. This corresponds to the second
case of Lemma \ref{lem-6-3}. Therefore the magnitude of the covariance term \eqref{eqn-6-11}
is bounded by $C/n_0^5$ for some constant $C$ in this case.

The number of such covariance terms is bounded by $Cn n_0^5$. The argument is the same as 
that in the previous case except for one difference. Given $i_1, k_1, l_1$, in this case,
we require one of $i_2,k_2,l_2$ to differ by $2$ or less  from some position $i^\ast$ such that
identical cards  occur at $i^\ast$ and at one of the positions in
$\nghd{i_1,i_1+1,k_1,l_1}$ in the deck $D_1$. Therefore, there are at most 
$Cn_0$ possibilities for one of $i_2,k_2,l_2$ and not just a constant number of possibilities
as in the previous case.

\item  Suppose $\nghd{i_1, i_1+1, k_1, l_1}\cap \nghd{i_2,i_2+1,k_2,l_2} = \phi$, and
$D_1(i^\ast_1) \neq D_1(i^\ast_2)$ for any $i^\ast_1 \in \nghd{i_1,i_1+1,k_1,l_1}$
and any $i^\ast_2\in\nghd{i_2,i_2+1,k_2,l_2}$. This corresponds to the third
case of Lemma \ref{lem-6-3}. Therefore the covariance term \eqref{eqn-6-11}
is $0$ in this case. 

\end{enumerate}

Consequently,
the total contribution of covariance terms of the form \eqref{eqn-6-11} to 
$\var(\mu Q)$ is bounded by $ C n$ for some positive constant $C$. Apart
from \eqref{eqn-6-11}, covariance terms can also be of the form
\begin{equation*}
\covar\Biggl(\chi_\pi(i_1,i_1+1), \frac{\chi_\pi(i_2,i_2+1,k,l)}{n_0^2}\Biggr),
\end{equation*}
with $i_1\in A$ and $i_2\in B$, or
\begin{equation*}
\covar\bigl(\chi_\pi(i_1,i_1+1), \chi_\pi(i_2,i_2+1)\bigr),
\end{equation*}
with $i_1\in A$ and $i_2\in A$. The proof that the total contribution of such
terms to $\var(\mu Q)$ is also bounded by $Cn$ is similar to and simpler than
the case that has already been dealt with. The bound on $\var(\mu Q)$ is stated
as a lemma below.

\begin{lem}
$\var(\mu Q) < Cn$ for some positive constant $C$.
\label{lem-6-4}
\end{lem}

\begin{thm}
Let $D_1$ be a deck of cards with $h$ different labels with each label occurring $n_0$ times.
Let $D_2$ be another deck with the same $n=hn_0$ cards in a different order. 
Let $\pi$ be a uniformly distributed permutation in $\Pi(D_1;D_2)$. Let the random
variable $W$ be the number of descents of $\pi$. Let $\mu = \E W$
and $\sigma^2 =\var(W)$. Assume $\sigma > 7^{2/3} (6 \mu)^{1/3}$.
Then 
\begin{equation*}
\Biggl| \P\Bigl(\frac{W-\mu}{\sigma}\leq x\Bigr) - \Phi(x)\Biggr|
\leq C \Bigl(\frac{1}{\sigma} + \frac{\mu}{\sigma^3} + \frac{\mu}{\sigma^4}
+ \frac{\sqrt{n}}{\sigma^2}\Bigr),
\end{equation*}
for some positive constant $C$ and $\Phi(x) = \frac{1}{\sqrt{2\pi}}\int_{-\infty}^x
\exp(-u^2/2) du$.
\label{thm-6-5}
\end{thm}
\begin{proof}
By Lemma \ref{lem-6-2}, $W^\ast$ has the $W$-size biased distribution. By construction
of $W^\ast$, $\abs{W^\ast-W}\leq 7$. By \cite[p. 477]{BillingsleyBook},
$\var\bigl(\E\bigl(W^\ast-W\bigl| W\bigr)\bigr)\leq \var(Q)$ with $Q$ defined by
\eqref{eqn-6-4}. The proof can be completed using Theorem \ref{thm-5-3}, if 
Lemma \ref{lem-6-4} is used to note that $\var(Q)\leq C n/\mu^2$.
\end{proof}

Consider the atypical case where $D_1= (1,2)^{n_0}$ and $D_2 = 1^{n_0}, 2^{n_0}$. Then every permutation in
$\Pi(D_1;D_2)$ has exactly $n_0-1$ descents, with each descent corresponding to
an occurrence of $2$ immediately before a $1$ in the deck $D_1$. In such a case
$\sigma^2 = 0$, but typically $\sigma^2$ will be of the order of $n$. As
$\mu < n$ always, in such cases, the normal approximation to $W$ given by Theorem 
\ref{thm-6-5} will have $O\bigl(n^{-1/2}\bigr)$ error. The computation of 
$\mu$ and 
$\sigma$
given $D_1$ and $D_2$ will be described presently.   

\vspace*{.2cm}
\noindent{\bf Calculating the mean and the variance of $W$.}
The expectation $\E W$ can be computed by summing $\E X_i$ over $1\leq i\leq n-1$.
If $D_1(i)=a$ and $D_1(i+1)=b$, the expectation $\E X_i$ is $0.5$
if $a=b$. If $a\neq b$, 
\begin{equation*}
\E X_i = \frac{R(a,b)}{N(a) N(b)},
\end{equation*}
where $N(a)$ and $N(b)$ are the number of cards with labels $a$ and
$b$ in the deck $D_2$, and $R(a,b)=\sum_{1\leq l_b < l_a\leq n}
\chi(l_a, l_b)$ with $\chi(l_a, l_b) = 1$ if $D_2(l_a) = a$
and $D_2(l_b) = b$ and $\chi(l_a, l_b)=0$ otherwise.

The variance $\var(W)$ can be obtained from $\E(X_i X_j)$ for
$1\leq i<j\leq n-1$. The computation of these joint expectations involves
many cases. First suppose that $j > i+1$. Denote $D_1(i)$, $D_1(i+1)$,
$D_1(j)$, and $D_1(j+1)$ by $a$, $b$, $c$, and $d$, respectively.
If $a=b$, then $\E(X_i X_j) = \E X_j/2$. Likewise if $c=d$,
$\E(X_i X_j) = \E X_i/2$. There are seven cases when $a\neq b$
and $c\neq d$. One of these  is when $a=d$ and $a,b,c$ are distinct.
In that case,
\begin{equation*}
\E(X_i X_j) = \frac{R(a,b) R(c,a) - R(c,a,b)}{N(a)N(b)N(c)(N(a)-1)},
\end{equation*}
where $R(c,a,b) = \sum_{1\leq l_c < l_a < l_b\leq n-1}\chi(l_c, l_a, l_b)$
with $\chi(l_c, l_a, l_c)=1$ if $D_2(l_c) = c$, $D_2(l_a)=a$ and
$D_2(l_b)=b$ but $\chi(l_c, l_a, l_c)=0$ otherwise. The other cases
are handled similarly. The $j=i+1$ case is also handled similarly.

\section{Shuffling cards for blackjack and bridge}
The rules for blackjack vary with the gambling house, but
the distinction between the four suits is always ignored. 
Most of the time all the face
cards are equivalent to cards with the number $10$, but there are some obscure situations
where $10$s, jacks, queens, and kings must all be considered distinct. 
We take a blackjack source deck to be any permutation of
the multiset $\{1^4,\ldots, 13^4\}$. We  consider two of these ---
$1^4,\ldots,13^4$ and $(1,\ldots,13)^4$ --- which are  notable
for their symmetry. These two will
be called {\it Blackjack1} and {\it Blackjack2}, respectively. 
Let $p_i$ be the transition
probability from one of these source decks to the $i$th possible ordering
of that deck after a certain number of riffle shuffles.
Ideally,
we would like all of these $p_i$ to be equal.

The situation for bridge is different. All the $52$ cards in the
source deck are distinct but there are only four players. Each player
must be dealt a random set of $13$ cards but the order in which
a player receives his cards is immaterial.

 Suppose the cards are dealt to players ${\mathcal N}$, ${\mathcal
 E}$, ${\mathcal S}$, and ${\mathcal W}$ in cyclic order, as is the
 common practice. Let the deck which needs to be dealt to the four
 players be $1,2,\ldots,52$. Let the sets $\mathcal{P}_{\mathcal N}$,
 $\mathcal{P}_{\mathcal E}$, $\mathcal{P}_{\mathcal S}$, and $\mathcal{P}_{\mathcal W}$ be a
partition
of $\{1,2,\ldots, 52\}$ with the cardinality of each set being $13$.
There are $52!/13!^4$ such partitions. Ideally, we would
like the probability that
${\mathcal N}$, ${\mathcal E}$, ${\mathcal S}$, and ${\mathcal W}$
receive cards in 
$\mathcal{P}_{\mathcal N}$, $\mathcal{P}_{\mathcal E}$, $\mathcal{P}_{\mathcal S}$,
 and $\mathcal{P}_{\mathcal W}$,
respectively, to be equal for all those partitions.

Suppose that the probability that ${\mathcal N}$, ${\mathcal E}$, ${\mathcal S}$, and ${\mathcal W}$
receive cards in 
$\mathcal{P}_{\mathcal N}$, $\mathcal{P}_{\mathcal E}$, $\mathcal{P}_{\mathcal S}$,
and $\mathcal{P}_{\mathcal W}$,
respectively, is equal to $p_\mathcal{P}$, when the deck 
$1,2,\ldots, 52$ is $a$-shuffled and then dealt
to ${\mathcal N}$, ${\mathcal E}$, ${\mathcal S}$, and ${\mathcal W}$
is cyclic order. Let the set of permutations of $\{1,2,\ldots,52\}$
that result in such a deal be $S_\mathcal{P}$. Now consider the deck $D$
with $52$ cards such that if $i$ belongs to 
$\mathcal{P}_{\mathcal N}$, $\mathcal{P}_{\mathcal E}$, $\mathcal{P}_{\mathcal S}$, 
or $\mathcal{P}_{\mathcal W}$, 
the $i$th card of $D$ is ${\mathcal N}$, ${\mathcal E}$, ${\mathcal S}$, 
or ${\mathcal W}$, respectively. Then the set of permutations $S_\mathcal{P}$
equals $\Pi(D; (\mathcal{NESW})^{13})$. Therefore, $p_\mathcal{P}$ is equal
to the probability that an $a$-shuffle of $D$ results in 
$(\mathcal{NESW})^{13}$. Ideally, we would like these transition 
probabilities to be equal for all $52!/13!^4$ possible decks $D$. 

In both situations, we have $N$ probabilities $p_i$, $1\leq i\leq N$,
which sum to $1$, and ideally we would like all of them to be $1/N$. 
Unfortunately,
 that
would require an infinite number of riffle shuffles. We have to determine how
close the probabilities are to the uniform distribution after a certain number
of riffle shuffles, and then decide how close is close enough.

 There are many ways to measure the closeness between probability
 distributions, but the notion of closeness must be chosen with 
 care. Some metrics such as the Euclidean distance or the $L^2$ norm
 are inappropriate as the following example shows. 
 Suppose one probability distribution always
 picks the first out of $N$ possibilities and another always picks the
 second. The Euclidean distance between these two distributions is
 $\sqrt{2}$. On the other hand, suppose a distribution picks one of the even
 numbered possibilities with equal likelihood and another
 distribution picks one of the odd possibilities with equal
 likelihood. The Euclidean distance between the third and fourth
 distributions is $2/\sqrt{N}+O(1/N)$, which incorrectly suggests that
 these two distributions are much closer to each other for large $N$.

\subsection{ A table of total variation distances}

The notion of closeness we use here is the {\it total variation distance}.
The total variation
distance from the probability distribution defined by the $p_i$ to the uniform
distribution is given by
\begin{equation}
\sum_{i=1}^N \Bigl(\frac{1}{N}-p_i\Bigr)^+,
\label{eqn-7-1}
\end{equation}
where $x^+ = x $ if $x\geq 0$ and $x^+=0$ if $x<0$. This total variation distance
always lies between $0$ and $1$. It also has a probabilistic meaning --- if $P_1(A)$
is the probability of a certain subset $A$ of the $N$  possibilities under the
distribution given by the $p_i$ and if $P_2(A)$ is the probability of that set
under the uniform distribution, then the total variation distance equals
the maximum of $\abs{P_1(A)-P_2(A)}$ over all possible $A$.
 
\begin{table}
\begin{center}
\scriptsize
\begin{tabular}{c|c|c|c|c|c|c|c|c|c|c}
 & 1 & 2 & 3 & 4 & 5 &  6 & 7 & 8 & 9 & 10\\ \hline
&&&&&&&&&&\\
BayerDiaconis &  {\bf 1} & {\bf 1} & {\bf 1} & {\bf 1} & 
{\bf .924} & {\bf .614} & {\bf .334} & {\bf .167} & {\bf .085} & {\bf .043}\\ 
&&&&&&&&&&\\
\hline
&&&&&&&&&&\\
Blackjack1 & {\bf 1} & {\bf 1} & {\bf 1} & {\bf .481} & 
{\bf .215} /.23 & {\bf .105} /.11 & {\bf .052} /.05 & {\bf .026} /.03 & {\bf .013} /.01 & {\bf .007} /.00\\ 
&&&&&&&&&&\\
\hline
&&&&&&&&&&\\
Blackjack2 & & & & &
.60 & .32 & .16 & .08 & .04 & .02 \\ 
&&&&&&&&&&\\
\hline
&&&&&&&&&&\\
Bridge1 & {\bf 1} & {\bf 1} & {\bf 1} & {\bf .990} &
{\bf .748}/.75 & {\bf .423}/.42 & {\bf .218}/.21 & {\bf .110}/.11 & {\bf .055}/.05 & {\bf .027}/.03\\ 
&&&&&&&&&&\\
\hline
&&&&&&&&&&\\
Bridge2 & & & .45? & .16 & 
.08 & .04 & .02 & .01 & .00 & .00\\ 
&&&&&&&&&&\\
\hline
&&&&&&&&&&\\
RedBlack1 & {\bf .580} & {\bf .360} & {\bf .208} & {\bf .105} & 
{\bf .052}/.05 & {\bf .026}/.03 & {\bf .013}/.01 & {\bf .007}/01 & {\bf .003}/.00 & {\bf .002}/.00\\ 
&&&&&&&&&&\\
\hline
&&&&&&&&&&\\
RedBlack2 & & & & &
.10 & .03 & .01 & .01 & .00 & .00\\ 
&&&&&&&&&&\\
\hline
&&&&&&&&&&\\
AliceBob1 & {\bf 1} & {\bf 1} & {\bf .999} & {\bf .725} &
{\bf .308} /.31 & {\bf .130}/.13 & {\bf .059}/.06 & {\bf .028}/.03 & {\bf .013}/.01 & {\bf .007}/.01\\ 
&&&&&&&&&&\\
\hline
&&&&&&&&&&\\
AliceBob2 & & & & &
.02 & .01 & .00 & .00 & .00 & .00
\end{tabular}
\normalsize
\end{center}
\caption[xyz]{Table of total variation distances after $1$ to $10$ riffle
shuffles for nine scenarios from BayerDiaconis to AliceBob2.  The
numbers in boldface have an error less that $.01$ with
a probability greater than $99.9996\%$. The numbers that
are not in boldface were determined using a less accurate and heuristic
method. Entries of the table with a number in boldface and in ordinary
type separated by a slash can be used to form an idea of the
reliability of the less accurate method.}
\label{table-1}
\end{table}


Table \ref{table-1} gives the total variation distances from the
uniform distribution for nine scenarios termed {\it BayerDiaconis},
{\it Blackjack1}, {\it Blackjack2}, {\it Bridge1}, {\it Bridge2},
{\it RedBlack1}, {\it RedBlack2}, {\it AliceBob1}, and {\it AliceBob2}.
That table allows the number of riffle shuffles
(\ie, $2$-shuffles) to vary from $1$ to
$10$.  In five of the scenarios, the source deck is fixed and the
target deck is allowed to vary. These
are {\it BayerDiaconis}, where the source deck is fixed
as $D_1 = 1,2,\ldots,52$; {\it Blackjack1} with
$D_1 = 1^{4}, 2^{4}, \ldots, {13}^4$; {\it Blackjack2} with
$D_2 = (1,2,\ldots,13)^4$; {\it RedBlack1} with
$D_1 = R^{26} B^{26}$; and {\it RedBlack2} with
$D_1 = (RB)^{26}$. In the other four scenarios the target deck is fixed.
These are {\it Bridge1}, where the target deck is fixed as
$D_2 = \mathcal{N}^{13} \mathcal{E}^{13} \mathcal{S}^{13} \mathcal{W}^{13}$; {\it Bridge2} with
$D_2 = (\mathcal{NESW})^{13}$; {\it AliceBob1} with $D_2 = A^{26} B^{26}$;
and {\it AliceBob2} with $D_2 = (AB)^{26}$.

How were the total variation distances shown in Table \ref{table-1}
determined? For {\it BayerDiaconis} there is a simple and elegant formula for
the total variation distance after a certain number of riffle shuffles
that was derived by Bayer and Diaconis \cite{BD1}. Some of the other numbers
require more extensive computation. Two numbers for {\it Bridge2}, with the number
of riffle shuffles being $3$ or $4$, were determined using more than
10000 hours of CPU time on the Teragrid computer network. 
From the point of view of card players,
those two numbers are perhaps the most interesting results of this paper.
First, we explain how the numbers in the rows of Table \ref{table-1},
other than the first {\it BayerDiaconis} row, were computed, and why they
can be trusted.

\subsection{Monte Carlo estimation of total variation distance}

The determination of the numbers in Table \ref{table-1}
is impeded by two problems. The first of these
is that the summation in \eqref{eqn-7-1} has so many terms that it is
impractical to use \eqref{eqn-7-1} to find the total variation distances
shown in Table \ref{table-1}.

 This first problem is
easy to overcome. Denote the sum in \eqref{eqn-7-1} by $S$. 
Let $X_1$ be a random variable that is equal to
$\bigl(1/N-p_i)^+$, the $i$th term in \eqref{eqn-7-1}, with probability
$1/N$ for $1\leq i\leq N$. Such a random variable can be easily generated
{\it if} we can determine the transition probabilities $p_i$ efficiently. 
We have $\E X_1 = S/N$ and $\var(X_1)\leq 1/N^2$. Let
$X_1, \ldots, X_k$ be  independent and identically distributed,
and let $Y_k = (X_1+\cdots+X_k)N/k$. Then 
\begin{equation}
\E Y_k = S.
\label{eqn-7-2}
\end{equation}
The theorem below tells us how good an estimate of $S$ can be
obtained from a single instance of the random variable $Y_k$.

\begin{thm}
For $\alpha> 0$ and $k>2$,
\[\P\Bigl(\abs{Y_k-S} \geq \frac{\alpha}{\sqrt{k}}\Bigr)
< \frac{4}{\alpha^4}.\] 
\label{thm-7-1}
\end{thm}
\begin{proof}
Consider the following calculation.
\begin{align*}
\E(\abs{Y_k-S}^4) &= \E\Biggl(\ABS{\Bigl(X_1-\frac{S}{N}\Bigr)
+\cdots + \Bigl(X_k-\frac{S}{N}\Bigr)}^4 \Biggr).\frac{N^4}{k^4}\\
&= \frac{N^4}{k^4}\sum_{i=1}^k \E\Bigl(X_i-\frac{S}{N}\Bigr)^4
+ \frac{6N^4}{k^4} \sum_{1\leq i < j\leq k}
 \E\Bigl(X_i-\frac{S}{N}\Bigr)^2\Bigl(X_j-\frac{S}{N}\Bigr)^2\\
&\leq (1/k^3+3/k^2) < 4/k^2
\end{align*}
The second equality above follows from $\E(X_i-S/N) = 0$ and from the
independence of $X_i$ and $X_j$ for $i\neq j$. To deduce the
first inequality in the last line above, note that $X_i$ has the
range $[0, 1/N]$ with $\E X_i = S/N$. The proof can be completed
using Markov's inequality \cite{BillingsleyBook}.
\end{proof}

In some special cases that include {\it Blackjack1}, {\it Bridge1},
{\it RedBlack1}, and {\it AliceBob1}, we have derived efficient polynomial time
algorithms for finding the descent polynomials \cite{CV1}. Those algorithms and
\eqref{eqn-3-2} can be used to determine the transition probabilities
efficiently. The numbers given for these four cases in Table 
\ref{table-1} in {\bf boldface}
were obtained using $k=10^7$ random permutations of the source deck or
the target deck, exact computation of the transition probabilities,
and Monte Carlo summation \eqref{eqn-7-2}.  Theorem \ref{thm-7-1} with
$k=10^7$ and $\alpha = \sqrt{10}$ implies that the boldface numbers
have errors less than $.001$ with probability greater than
$96\%$. Theorem \ref{thm-7-1} with $k=10^7$ and $\alpha = 10\sqrt{10}$
implies that the boldface numbers have errors less than $.01$ with
probability greater than $99.9996\%$.

The more significant
problem 
is that there may be no efficient way to determine the transition
probabilities $p_i$. Indeed there is very probably no efficient method for
determining these transition probabilities in general, as proved in Section 4.

\subsection{Approximation of the descent polynomial}

In the other four cases --- {\it Blackjack2}, {\it Bridge2}, {\it RedBlack2},
and {\it AliceBob2} --- we turn to the Monte Carlo method once again to
approximate the descent polynomials. Suppose we are given decks $D_1$
and $D_2$ and it is required to approximate the descent polynomial of
permutations in $\Pi(D_1;D_2)$. If the decks $D_1$ and $D_2$ have
$n_c$ cards with label $c$ for $1\leq c \leq h$, then the total number 
of permutations in $\Pi(D_1;D_2)$ is $n_1!\ldots n_h!$. We generate
$l$ random permutations $\pi_1,\ldots, \pi_l$ from this collection and form
the polynomial
\begin{equation}
P = \sum_{i=1}^l x^{\des(\pi_i)} = \gamma_0+\gamma_1 x+\cdots+\gamma_{n-1}
x^{n-1}.
\label{eqn-7-3}
\end{equation}
The coefficient $\gamma_d$ counts the number of random permutations 
with $d$ descents.
The approximation to the descent polynomial is taken to be
$\frac{n_1!\ldots n_h!}{l} P$; in other words, the polynomial $P$
given by \eqref{eqn-7-3} is normalized to get an approximation to the descent 
polynomial. 

Once the polynomial $P$ defined by \eqref{eqn-7-3} is formed, it carries within
itself an estimate of the accuracy of the coefficients of the descent polynomial
obtained by normalizing $P$, as we will explain.

Suppose that $X$ is a Bernoulli random variable with $\P(X=1) =
p$ and that $p$ is unknown. 
Suppose that $X_1, \ldots X_l$ are independent with the same
distribution as that of $X$, and that in one experiment $m$ out of these
$l$ random variables equal $1$. We can estimate $p \approx m/l$, but
how accurate is this estimate?
Let $Y = (X_1+\ldots + X_l)/l$. Then $\E Y = p$ and $\var(Y) = p(1-p)/l$.
Therefore the fluctuations of $Y$ about its mean are of the order $\sqrt{p(1-p)/l}$.
If we use a single instance of $Y$ to estimate its mean, which is $p$, then we
expect an absolute error of about $\sqrt{p(1-p)/l}$ and a relative error
of about $\sqrt{(1-p)/l p}$. If we substitute $p = m/l$, we find that
the relative error will be about $\sqrt{(l-m)/l m}$. If $p$ is very small,
then $m << l$ and we may expect a relative error of about $1/\sqrt{m}$.

If we define $X(\pi) = 1$ if $\des(\pi) = d$ and $X(\pi) = 0$
otherwise, where $\pi$ is a uniformly distributed permutation in
$\Pi(D_1;D_2)$, then $\gamma_d$ defined by
\eqref{eqn-7-3} is $X_1+ \cdots+ X_l$, where $X_i=X(\pi_i)$ are independent
with the same distribution as that of $X(\pi)$. 
Then by the argument of the preceding paragraph, the relative
error in the estimate of the $x^d$ coefficient of the descent
polynomial will be about $1/\sqrt{\gamma_d}$.

To illustrate this estimate in practice, we take the target deck to be
$D_2 = (\mathcal{NESW})^{13}$, which is fixed for {\it Bridge2}, and the source
deck $D_1$ to be
\begin{quote}
\footnotesize
\cal{NSEENNWEWSSWESWNNNEESSSSSESWWNNSENWSEWSWWWEENEWNNNWE}.
\normalsize
\end{quote}
We computed $P$, which is defined by \eqref{eqn-7-3}, with $l=10^9$,
and got the coefficients of $x^{14}$, $x^{15}$, and $x^{16}$ to be
$17$, $397$, and $4560$, respectively. If the descent polynomial of
permutations from $D_1$ to $D_2$ is approximated as $13!^4 P/10^9$, we
expect the relative errors in the coefficients of $x^{14}$, $x^{15}$,
and $x^{16}$ to be about $25\%$, $5 \%$, and $1.6\%$, respectively.
When we computed $P$ with $l = 10^{11}$, we got the coefficients
of $x^{14}$, $x^{15}$, and $x^{16}$ to be $2334$, $38418$, 
and $468359$, respectively. These numbers are $100$ times the
counts for $l=10^{9}$, if allowance is made for the relative errors
when the counts are scaled and interpreted as coefficients
of the descent polynomial. We conclude that  descent polynomials can be
approximated using \eqref{eqn-7-3} and that the accuracy
of the approximation can be gauged by looking at the coefficients of $P$.

By  \eqref{eqn-2-1}, the probability that an
$a$-shuffle leads to a permutation with $d$ descents is equal to
\[\frac{1}{a^n}\binom{a+n-d-1}{n} \Euler{n}{d},\]
where $\euler{n}{d}$ is the Eulerian number that counts
the number of permutations of $1,2,\ldots,n$
with $d$ descents. The Eulerian numbers can be calculated using
simple recurrences \cite{GKPBook}.

The probability that an $a$-shuffle
with $a=32$ has $16$ or more descents is more than $0.95$.
Coefficients of terms from
$x^{16}$ to $x^{36}$ of descent polynomials that correspond to the scenarios
from {\it Blackjack1} to {\it AliceBob2} in Table \ref{table-1} can
be approximated well using \eqref{eqn-7-3} with $l=10^7$.
The numbers reported in Table \ref{table-1} for these scenarios, with
the number of riffle shuffles varying from $5$ to $10$,
that are not in boldface
 were computed
by approximating descent polynomials in this manner. The transition
probabilities were computed using these approximate descent polynomials
and \eqref{eqn-3-2}, and the total variation distances were computed using
\eqref{eqn-7-2} with $k=1000$. The total variation distances computed this
way are reported with two digits after the decimal point. They compare
well with more accurate estimates of the total variation distances,
which are found in boldface in Table \ref{table-1},
for {\it Blackjack1}, {\it Bridge1}, {\it RedBlack1}, and {\it AliceBob1}.

We still need to explain the computation of the total variation
distances for {\it Bridge2} when the number of riffle shuffles is $3$
or $4$.  Four riffle shuffles are equivalent to an $a$-shuffle with
$a=16$, and it is necessary to accurately compute the coefficients of
terms from $x^{12}$ to $x^{16}$ of the descent polynomials to find the
total variation distance for {\it Bridge2} after four riffle
shuffles. Obtaining an accurate estimate for the coefficient of
$x^{12}$ using \eqref{eqn-7-3} would require an $l$ that is beyond the
reach of today's computers. We used \eqref{eqn-7-3} with $l=10^{10}$,
and with this $l$ the coefficient of $x^{15}$ is approximated with
a relative error that is less than $5\%$ with a probability greater
than $95\%$. We took the $\log$s of the coefficients of the terms from
$x^{15}$ to $x^{22}$ and computed degree $4$ polynomials that were
least squares fits to these numbers. These polynomials were nearly
quadratic in accord with Theorem \ref{thm-6-5}. We got the
coefficients of the $x^{12}$, $x^{13}$, and $x^{14}$ terms by
extrapolation. 
We feel sure that the extrapolated coefficients had relative errors
smaller than $10\%$.

For {\it Bridge2} and three riffle shuffles, we got the coefficients
of the $x^7$ terms in the descent polynomials using the coefficients
of terms from $x^{15}$ to $x^{26}$ and polynomial fits of degree $6$.
We are less sure that the estimated total variation distance 
for this case is accurate.

\section{Three open problems}

\noindent {\bf 1.} The first of the three open problems that we mention in this section
is about a more general model of riffle shuffles. In the model of
riffle shuffles we have used so far a random riffle shuffle of
a deck of $n$ cards is obtained by first generating a random 
sequence of $n$ numbers. The numbers in the sequence are independent
of each other and each number is either $1$ or $2$ with probability $1/2$.
The model can be changed
by requiring the first number in the sequence to be either $1$ or $2$
with probability $1/2$. Every later number in the sequence equals the
preceding number with probability $1-p$ and it is of the opposite
kind with probability  $p$. This model is described
by Aldous \cite{Aldous1} and Diaconis \cite{Diaconis2}.

When $p=1/2$, we get back the GSR-model which was described in Section 2.
If $p>1/2$, then the riffle shuffles are neat --- if a card is dropped
from one hand the next card is likelier to be dropped from the other hand.
If $p<1/2$, the shuffling is clumsy. The problem is to determine how
the mixing times depend upon $p$.

Any sequence of $1$s and $2$s where every $1$ occurs before every $2$
corresponds to the identity permutation. All other sequences correspond
to distinct permutations. The probability of one of these permutations under
this model will be  $0.5*p^k(1-p)^{n-1-k}$ if there are $k$ places
where a $2$ follows a $1$ or a $1$ follows a $2$ in the
corresponding sequence. The probability of the identity permutation
is also easy to determine. Suppose we want to know the probability that
a given permutation $\pi$ can be obtained as the composition of
$m$ riffle shuffles none of which is the identity. This number will
be a polynomial in $p$ of the form
\begin{equation*}
\sum_{k=0}^{m(n-1)} c_k p^k (1-p)^{m(n-1)-k}
\end{equation*}
with the coefficients $c_k$ depending upon $\pi$ and $m$. A good place
to start might be by asking if the coefficients
$c_k$ can be determined in polynomial time.

The value of this polynomial when $p=1/2$ can be deduced from the
work of Bayer and Diaconis \cite{BD1}. When $p=1$, each riffle shuffle
is either a perfect in-shuffle or a perfect out-shuffle with 
probability $1/2$. Much  information about this case can be found
in the work of Diaconis, Graham, and Kantor \cite{DGK1}.

\vspace*{0.2cm}
\noindent {\bf 2.} 
The second problem is purely combinatorial and is related to the
riffle shuffles of the deck $(1,2)^n$ \cite{CV1}. Consider all
sequences of length $2n$ with $n$ $1$s and $n$ $2$s. Define any two
sequences $D_1 = \alpha \beta \gamma$ and $D_2 = \alpha \beta^\ast
\gamma$ to be $R$-related if $\beta$ has the same number of $1$s as
$2$s and if $\beta^\ast$ is obtained from $\beta$ by changing $1$s to
$2$s and $2$s to $1$s. This is an equivalence relation. We conjecture
that the number of equivalence classes is $(n+3) 2^{n-2}$.  If the
equivalence relation is modified by requiring $\alpha$ and
$\gamma$ to be sequences of the same length, the number of equivalence
classes is the Catalan number $\frac{1}{n+2}\binom{2n+2}{n+1}$ \cite{CV1}.

\vspace{0.2cm}
\noindent {\bf 3.}
The third problem too is purely combinatorial. Consider all
permutations $\pi$ of the numbers $1,\ldots, nh$ such that 
$\pi(i)\equiv i \mod h$ for $1\leq i \leq nh$. The problem is to
derive a recurrence for the number of these permutations that 
have exactly $d$ descents. If $h=1$, the familiar recurrences for
the Eulerian numbers (see \cite{GKPBook}) solve this problem. A solution
of this problem will make it possible to compute the transition
probability from the deck $(1,2,\ldots,h)^n$ to itself under an
$a$-shuffle.

\section{Acknowledgements} The authors thank P. Diaconis, S. Fomin,
J. Fulman, C. Mulcahy, and J. Stipins for helpful discussions. In
addition, the authors thank Prof. Diaconis for his work, which got
them interested in this problem, and for telling them about Stein's
method.


\bibliography{references}
\bibliographystyle{plain}
\end{document}